\documentclass{article}
\pdfoutput=1

\usepackage{arxiv}

\usepackage[utf8]{inputenc} 
\usepackage[T1]{fontenc}    
\usepackage{hyperref}       
\usepackage{url}            
\usepackage{booktabs}       
\usepackage{amsfonts}       
\usepackage{nicefrac}       
\usepackage{microtype}      

\usepackage{lipsum}         
\usepackage{graphicx}
\usepackage{natbib}
\usepackage{doi}

\usepackage{mathtools,amsthm}
\usepackage{algorithm}
\usepackage{algpseudocode}
\usepackage{color}
\allowdisplaybreaks[4]
\usepackage{cleveref}       
\usepackage{etoolbox}

\usepackage{bbm}
\usepackage{amsthm}
\usepackage{natbib}
\usepackage{mathtools}
\usepackage{algorithm}
\usepackage{algpseudocodex}
\usepackage{cleveref}  
\usepackage{makecell}
\usepackage{nicematrix}
\usepackage{enumitem}
\usepackage{amssymb}
\usepackage{tablefootnote}
\usepackage{soul}
\usepackage{scalefnt}

\hypersetup{
	colorlinks=true,
	linkcolor=red!75!black,
	citecolor=blue!50!black
}

\newtheorem{theorem}{Theorem}
\newtheorem{lemma}{Lemma}

\newtheorem{assumption}{Assumption}

\newtheorem{definition}{Definition}
\crefname{assumption}{assumption}{assumptions}

\DeclarePairedDelimiter{\abs}{\lvert}{\rvert}
\DeclarePairedDelimiter{\norm}{\|}{\|}
\DeclarePairedDelimiter{\sqn}{\|}{\|^2}

\DeclarePairedDelimiter{\floor}{\lfloor}{\rfloor}

\def\<#1,#2>{\langle #1,#2\rangle}

\DeclareMathOperator{\dom}{dom}

\DeclareMathOperator{\spanset}{span}

\newcommand{\R}{\mathbb{R}}

\newcommand{\ones}{\mathbf{1}}

\newcommand{\lmax}{\lambda_{\max}}

\newcommand{\lminp}{\lambda_{\min}^+}

\newcommand{\cE}{\mathcal{E}}

\newcommand{\cG}{\mathcal{G}}

\newcommand{\cK}{\mathcal{K}}
\newcommand{\cL}{\mathcal{L}}
\newcommand{\cM}{\mathcal{M}}

\newcommand{\cO}{\mathcal{O}}

\newcommand{\cV}{\mathcal{V}}

\newcommand{\mI}{\mathbf{I}}

\newcommand{\mP}{\mathbf{P}}

\newcommand{\mW}{\mathbf{W}}

\newcounter{annotatecount}
\newcounter{annotateidx}
\newcounter{annotatejdx}
\newcounter{annotatelabelcount}
\newcommand{\atran}[2]{\stepcounter{annotatecount}\overset{(\alph{annotatecount})}{#1}\csgdef{annotatedescription\theannotatecount}{#2}}
\newcommand{\aeq}[1]{\atran{=}{#1}}
\newcommand{\aleq}[1]{\atran{\leq}{#1}}

\newcommand{\asubset}[1]{\atran{\subset}{#1}}

\newcommand{\annotateinitused}{\setcounter{annotateidx}{0}\whileboolexpr{test{\ifnumless{\theannotateidx}{\theannotatecount}}}{\stepcounter{annotateidx}\csgdef{aused\theannotateidx}{0}}}
\newcommand{\annotategetlabels}{\setcounter{annotatejdx}{0}\setcounter{annotatelabelcount}{0}\whileboolexpr{test{\ifnumless{\theannotatejdx}{\theannotatecount}}}{\stepcounter{annotatejdx}\ifcsequal{annotatedescription\theannotateidx}{annotatedescription\theannotatejdx}{\csgdef{aused\theannotatejdx}{1}\stepcounter{annotatelabelcount}\csedef{annotatelabel\theannotatelabelcount}{(\alph{annotatejdx})}}{}}}
\newcommand{\annotateprintlabels}{\setcounter{annotatejdx}{0}\whileboolexpr{test{\ifnumless{\theannotatejdx}{\theannotatelabelcount}}}{\stepcounter{annotatejdx}\ifnumequal{\theannotatejdx}{\theannotatelabelcount}{\ifnumequal{\theannotatejdx}{1}{}{~and~}}{}\csuse{annotatelabel\theannotatejdx}\ifnumless{\theannotatejdx}{\theannotatelabelcount}{\ifnumless{\theannotatejdx+1}{\theannotatelabelcount}{,~}{}}{}}}
\newcommand{\annotate}{\annotateinitused\setcounter{annotateidx}{0}\whileboolexpr{test{\ifnumless{\theannotateidx}{\theannotatecount}}}{\stepcounter{annotateidx}\ifcsstring{aused\theannotateidx}{0}{\ifnumequal{\theannotateidx}{1}{}{;~}\annotategetlabels\annotateprintlabels~\csuse{annotatedescription\theannotateidx}}{}}\setcounter{annotatecount}{0}}

\newcommand{\sX}{{\R^{d}}}
\newcommand{\rng}[2]{\{#1,\ldots,#2\}}

\newcommand{\ol}[1]{\overline{#1}}
\newcommand{\uln}[1]{\underline{#1}}

\newcommand{\mind}[1]{\hspace{#1}&\hspace{-#1}}

\newcommand{\alertR}[1]{{\color{red}#1}}
\newcommand{\alertB}[1]{{\color{blue}#1}}

\newcommand{\OA}{$^\dagger$}
\newcommand{\LB}{$^*$}

\newcommand{\mem}[1]{\cM_{#1}}
\newcommand{\memloc}[1]{\cM_{#1}^{\text{sub}}}
\newcommand{\memcom}[1]{\cM_{#1}^{\text{com}}}
\newcommand{\tloc}{\tau_{\text{sub}}}
\newcommand{\tcom}{\tau_{\text{com}}}
\newcommand{\texe}{\tau}
\newcommand{\subgrad}{\hat{\nabla}}
\newcommand{\basis}[2]{{\bf e}_{#2}^{#1}}

\newcommand{\hlcom}[1]{\colorbox{green!35}{$\displaystyle #1$}}
\newcommand{\hlloc}[1]{\colorbox{yellow!60}{$\displaystyle #1$}}
\newcommand{\algname}[1]{{\sf\scalefont{0.92}{#1}}}

\title{Lower Bounds and Optimal Algorithms for Non-Smooth Convex Decentralized Optimization over Time-Varying Networks}

\date{}

\newif\ifuniqueAffiliation
\uniqueAffiliationtrue

\ifuniqueAffiliation 
\author{
	Dmitry Kovalev \\
	Yandex Research \\
	\texttt{dakovalev1@gmail.com}
	\And
	Ekaterina Borodich\\
	MIPT\thanks{Moscow Institute of Physics and Technology}\\
	\texttt{borodich.ed@phystech.edu}
	\AND
	Alexander Gasnikov\\
	Innopolis University, MIPT, ISP RAS\thanks{Institute for System Programming of the Russian Academy of Sciences}\\
	\texttt{gasnikov@yandex.ru}
	\And
	Dmitrii Feoktistov\\
	MSU\thanks{Moscow State University}\\
	\texttt{feoktistovdd@my.msu.ru}
}
\else
\usepackage{authblk}

\author[1]{%
	Dmitry Kovalev\thanks{\texttt{dakovalev1@gmail.com}}
}
\author[2]{
	Ekaterina Borodich\thanks{\texttt{borodich.ed@phystech.edu}}
}
\author[2]{
	Alexander Gasnikov\thanks{\texttt{gasnikov@yandex.ru}}
}
\author[3]{
	Dmitrii Feoktistov\thanks{\texttt{feoktistovdd@my.msu.ru}}
}
\affil[1]{Yandex Research}
\affil[2]{Moscow Institute of Physics and Technology}
\affil[3]{Moscow State University}
\fi
\renewcommand{\headeright}{}
\renewcommand{\undertitle}{}
\renewcommand{\shorttitle}{Optimal Algorithms for Non-Smooth Convex Decentralized Optimization over Time-Varying Networks}

\hypersetup{
pdftitle={Lower Bounds and Optimal Algorithms for Non-Smooth Convex Decentralized Optimization over Time-Varying Networks},
pdfsubject={math.OC},
pdfauthor={Dmitry Kovalev, Ekaterina Borodich, Alexander Gasnikov, Dmitrii Feoktistov},
pdfkeywords={decentralized optimization, convex optimization, distributed optimization, non-smooth optimization, time-varying networks},
}

\begin{document}
\maketitle

\begin{abstract}
	We consider the task of minimizing the sum of convex functions stored in a decentralized manner across the nodes of a communication network. This problem is relatively well-studied in the scenario when the objective functions are smooth, or the links of the network are fixed in time, or both. In particular, lower bounds on the number of decentralized communications and (sub)gradient computations required to solve the problem have been established, along with matching optimal algorithms. However, the remaining and most challenging setting of non-smooth decentralized optimization over time-varying networks is largely underexplored, as neither lower bounds nor optimal algorithms are known in the literature. We resolve this fundamental gap with the following contributions: (i) we establish the first lower bounds on the communication and subgradient computation complexities of solving non-smooth convex decentralized optimization problems over time-varying networks; (ii) we develop the first optimal algorithm that matches these lower bounds and offers substantially improved theoretical performance compared to the existing state of the art.
\end{abstract}


\section{Introduction}\label{sec:intro}

In this paper, we study the decentralized optimization problem. Specifically, given a set of $n$ compute nodes connected through a communication network, our goal is to solve the following finite-sum optimization problem with quadratic regularization:
\begin{equation}\label{eq:main}
	\min_{x \in \sX} \left[p(x) = \frac{1}{n}\sum_{i=1}^{n}f_i(x) + \frac{r}{2}\sqn{x}\right],
\end{equation}
where $r \geq 0$ is a regularization parameter, and each function $f_i(x)\colon \sX \to \R$ is stored on the corresponding node $i \in \rng{1}{n}$. Each node $i$ can perform computations based on its local state and data, and can directly communicate with other nodes through the links in the communication network.

Decentralized optimization problems find applications in a wide variety of fields. These include network resource allocation \citep{beck20141}, distributed model predictive control \citep{giselsson2013accelerated}, power system control \citep{gan2012optimal}, distributed spectrum sensing \citep{bazerque2009distributed}, and optimization in sensor networks \citep{rabbat2004distributed}. In addition, such problems cover the supervised training of machine learning models through empirical risk minimization, thus attracting significant interest from the machine learning community \citep{lian2017can,ryabinin2021moshpit,ryabinin2020towards}.

\subsection{Time-varying Networks}

In our paper, we focus on the setting in which the links in the communication network are allowed to change over time. Such time-varying networks \citep{zadeh1961time,kolar2010estimating} hold significant relevance to many practical applications. For instance, in sensor networks, changes in the link structure can be caused by the motion of sensors and disturbances in the wireless signal connecting pairs of sensors. Similarly, in distributed machine learning, connections between compute nodes can intermittently appear and disappear due to network unreliability \citep{ryabinin2020towards}. Lastly, we anticipate that the time-varying setting will be supported by future-generation federated learning systems \citep{konecny2016federated,mcmahan2017communication}, where communication between pairs of mobile devices or between mobile devices and servers will be affected by their physical proximity, which naturally changes over time.

\subsection{Convex Setting}\label{sec:intro_cvx}

In this work, we consider the decentralized optimization problem in the case when the objective function is convex (or strongly convex). At first glance, it may seem that the convexity assumption is restrictive and should not be considered. However, as we will see further, even in this fundamental setting, the existing algorithmic developments are limited and have significant gaps that need to be closed. Moreover, considering the convex optimization setting offers important benefits compared to general non-convex functions. One such benefit is that convex optimization often serves as a source of inspiration for the development of algorithms that turn out to be highly effective in solving practical problems, even non-convex ones.

For example, state-of-the-art optimization algorithms such as \algname{Adam} \citep{kingma2014adam} and \algname{RMSProp} \citep{hinton2012neural} employ the momentum trick, which is observed to be efficient for numerous tasks, including the training of deep neural networks. However, from the perspective of non-convex optimization theory, momentum is useless because, for non-convex problems, the iteration complexity of the standard gradient method cannot be improved \citep{carmon2020lower}. On the other hand, it was theoretically proven that momentum substantially boosts the convergence speed of the gradient method when applied to convex functions \citep{nesterov1983method}. In other words, convex optimization theory suggests that the momentum trick should be used, while non-convex theory suggests that it should not, and the former aligns much more closely with practical observations. A similar situation can be seen with other state-of-the-art optimization methods, including distributed local gradient methods \citep{mishchenko2022proxskip,sadiev2022communication,karimireddy2020scaffold}, adaptive gradient methods \citep{duchi2011adaptive}, etc.
Such inconsistency between non-convex theoretical convergence guarantees for optimization algorithms and their actual performance in practice can be attributed to the fact that the class of non-convex functions is far too broad. This is why many optimization research papers try to narrow down this class by considering additional assumptions such as Polyak-{\L}ojasiewicz condition \citep{karimi2016linear}, bounded non-convexity \citep{carmon2018accelerated,allen2018natasha}, quasi-strong convexity \citep{necoara2019linear}, etc. However, these assumptions can be seen as relaxations of the standard convexity property. Therefore, we naturally opt to focus on the convex decentralized optimization problem, leaving potential generalizations for future work.

\subsection{Related Work and Main Contributions}\label{sec:intro_contrib}

Decentralized optimization has been attracting a lot of attention for more than a decade. Plenty of algorithms have been developed, including \algname{EXTRA} \citep{shi2015extra}, \algname{DIGing} \citep{nedic2017achieving}, \algname{SONATA} \citep{scutari2019distributed}, \algname{NIDS} \citep{li2019decentralized}, \algname{APM-C} \citep{li2018sharp,rogozin2021towards}, and many others. In recent years, the focus of the research community has shifted towards the more complex task of finding, in some sense, the best possible algorithms for solving decentralized optimization problems \citep{scaman2017optimal,scaman2018optimal,lan2020communication,kovalev2020optimal,kovalev2021adom,kovalev2021lower,kovalev2022optimal,hendrikx2021optimal,li2022variance,li2021accelerated,metelev2024decentralized}. This task consists of finding a lower bound on the complexity\footnote{By complexity, we mean, depending on the context, the number of subgradient computations or decentralized communications required to solve the problem.} of solving a given subclass of decentralized problems and finding an algorithm whose complexity matches this lower bound. Such algorithms are called optimal because their complexity cannot be improved for a given problem class due to the established lower bounds.

We discuss the four main classes of decentralized optimization problems that cover smooth\footnote{A function is called smooth if it is continuously differentiable and has a Lipschitz-continuous gradient.} and non-smooth objective functions, and fixed and time-varying communication networks. We reference the existing state-of-the-art research papers that collectively solve the task of finding optimal algorithms for these classes. These papers are summarized in \Cref{table:main}. In the case of smooth and strongly convex objective functions and fixed communication networks, \citet{scaman2017optimal} established the lower bounds on the number of communication rounds and the number of local gradient computations required to find the solution. These lower bounds were matched by \algname{OPAPC} algorithm of \citet{kovalev2020optimal}. In the case of smooth and strongly convex problems over time-varying networks, lower complexity bounds were provided by \citet{kovalev2021lower}, and two optimal algorithms were developed: \algname{ADOM+} \citep{kovalev2021lower} and \algname{AccGT} \citep{li2021accelerated}. In the case of non-smooth convex problems over fixed networks, lower bounds were established by \citet{scaman2018optimal}, and two optimal algorithms were proposed: \algname{DCS} \citep{lan2020communication} and \algname{MSPD} \citep{scaman2018optimal}.

Our paper primarily focuses on the remaining and most challenging setting of non-smooth convex decentralized optimization problems over time-varying networks. Only a few algorithms have been developed for this setting, including the distributed subgradient method (\algname{D-SubGD}) by \citet{nedic2009distributed}, the subgradient-push method (\algname{SubGD-Push}) by \citet{nedic2014distributed}, and \algname{ZOSADOM} by \citet{lobanov2023non}. Moreover, to the best of our knowledge, neither lower complexity bounds nor optimal algorithms have been proposed in this setting. Consequently, in this work, we close this significant gap with the following key contributions:
\begin{enumerate}
	\item[\bf (i)] We establish the first lower bounds on the number of decentralized communications and local subgradient computations required to solve problem~\eqref{eq:main} in the non-smooth convex setting over time-varying networks,
	\item[\bf (ii)] We show that our lower bounds are tight by developing the first optimal algorithm that matches these lower bounds. The proposed algorithm has state-of-the-art theoretical communication complexity, which outclasses the existing methods described in the literature.
\end{enumerate}

\begin{table}[t]
	\centering
	\caption{Summary of the existing state-of-the-art results in decentralized convex optimization. Multiple paper references are provided for each problem setting: papers marked with \LB~provide lower complexity bounds, and papers marked with \OA~provide optimal algorithms that match the corresponding lower bounds.
	}
	\label{table:main}
	\begin{NiceTabular}{|c|c|c|}
		\CodeBefore
		\cellcolor{cyan!12}{3-3,3-3}
		\Body
		\toprule
		                   &
		\bf Smooth Setting &
		\bf Non-Smooth Setting
		\\\midrule
		\makecell{
			\bf Fixed
		\\\bf Networks
		}
		                   &
		\makecell{
			\citet{kovalev2020optimal}\OA
		\\\citet{scaman2017optimal}\LB
		}                  &
		\makecell{
			\citet{lan2020communication}\OA
		\\\citet{scaman2018optimal}\OA\LB
		}
		\\\midrule
		\makecell{
			\bf Time-Varying
		\\\bf Networks
		}
		                   &
		\makecell{
			\citet{kovalev2021lower}\OA\LB
		\\\citet{li2021accelerated}\OA
		}                  &
		\makecell{
			\bf\Cref{alg} (this paper)\OA
		\\\bf\Cref{thm:lb_sc,thm:lb} (this paper)\LB
		}
		\\\bottomrule
	\end{NiceTabular}
\end{table}

\section{Notation and Assumptions}\label{sec:ass}

In this paper, we are going to use the following notations:
$\otimes$ denotes the Kronecker matrix product,
$\mI_p$ denotes a $p\times p$ identity matrix,
$\ones_{p} = (1,\ldots,1) \in \R^p$,
$\basis{p}{j} \in \R^p$ for $j \in \rng{1}{p}$ denotes the $j$-th unit basis vector,
where $p \in \{1,2,\ldots\}$. In addition, $\norm{\cdot}$ denotes the standard Euclidean norm of a vector, and $\<\cdot,\cdot>$ denotes the standard scalar product of two vectors.

\subsection{Objective Function}

Further, we describe the assumptions that we impose on problem~\ref{eq:main}. As discussed in \Cref{sec:intro_cvx}, we assume the convexity of the objective function in problem~\eqref{eq:main}. In particular, we assume that functions $f_1(x),\ldots,f_n(x)$ are convex, which is formally described in \Cref{ass:cvx}.
\begin{assumption}\label{ass:cvx}
	Each function $f_i(x)$ is convex. That is, for all $x',x \in \sX$ and $\tau \in [0,1]$, the following inequality holds:
	\begin{equation}
		f_i(\tau x + (1-\tau) x') \leq \tau f_i(x) + (1-\tau)f_i(x').
	\end{equation}
\end{assumption}
In addition, we assume that the objective functions $f_1(x),\ldots,f_n(x)$ are Lipschitz continuous, which is formalized in \Cref{ass:lip}. This property is widely used in the theoretical analysis of non-smooth optimization algorithms, such as the subgradient method \citep{nesterov2013introductory}, dual extrapolation method \citep{nesterov2009primal}, etc.
\begin{assumption}\label{ass:lip}
	Each function $f_i(x)$ is $M$-Lipschitz continuous for $M \geq 0$. That is, for all $x', x \in \sX$, the following inequality holds:
	\begin{equation}
		\abs{f_i(x) - f_i(x')} \leq M \norm{x - x'}.
	\end{equation}
\end{assumption}
We also need the following \Cref{ass:sol}, which ensures the existence of a solution to problem~\eqref{eq:main}. Note that in the strongly convex case ($r > 0$), the solution always exists and is unique. However, in the convex case ($r=0$), we need to explicitly assume the existence of a solution.
\begin{assumption}\label{ass:sol}
	There exists a solution $x^* \in \sX$ to problem~\eqref{eq:main} and a distance $R > 0$ such that
	$
		\norm{x^*}\leq R.
	$
\end{assumption}

\subsection{Decentralized Communication}

Next, we formally describe the decentralized communication setting. The communication network is typically represented by a graph $\cG(\cV,\cE)$, where $\cV = \rng{1}{n}$ is the set of compute nodes and $\cE \subset \cV \times \cV$ is the set of links in the network. As mentioned earlier, we allow the communication links to change over time. Thus, we introduce the continuous time parameter $\texe \geq 0$ and a set-valued function $\cE(\texe)\colon \R_{+} \to 2^{\cV \times \cV}$, which represents the time-varying set of edges.\footnote{By $2^{\cV \times \cV} = \{\cE : \cE \subset \cV \times \cV\}$ we denote the set of all subsets of $\cV \times \cV$.} Our time-varying network is then denoted as $\cG(\texe) = (\cV,\cE(\texe))$.

Decentralized communication is typically represented via a matrix-vector multiplication with the so-called gossip matrix associated with the communication network \citep{scaman2017optimal,kovalev2021lower}. In the time-varying setting, we represent the gossip matrix by a matrix-valued function $\mW(\texe) \colon \R_+ \to \R^{n\times n}$, which satisfies the following \Cref{ass:gossip}.
\begin{assumption}\label{ass:gossip}
	For all $\texe \geq 0$, the gossip matrix $\mW(\texe) \in \R^{n\times n}$ associated with the time-varying communication network $\cG(\cV, \cE(\texe))$ satisfies the following properties:
	\begin{enumerate}
		\item[\bf (i)] $\mW(\texe)_{ij} = 0$\;\;if\;\;$i\neq j$ and $(j,i) \notin \cE(\texe)$,
		\item[\bf (ii)] $\mW(\texe) \ones_n = 0$\;\;and\;\;$\mW(\texe)^\top \ones_n = 0$.
	\end{enumerate}
\end{assumption}
We also define the so-called condition number of the network $\chi \geq 1$, which indicates how well the network $\cG(\texe)$ is connected \citep{scaman2017optimal,kovalev2021lower}. In particular, the communication complexity of most decentralized optimization algorithms depends on $\chi$. \Cref{ass:chi} provides the formal definition of this quantity.
\begin{assumption}\label{ass:chi}
	There exists a constant $\chi \geq 1$ such that the following inequality holds for all $\texe \geq 0$:
	\begin{equation}
		\sqn{\mW(\texe) x  - x} \leq \left(1 - {1}/{\chi}\right)\sqn{x}
		\;\;
		\text{for all}
		\;\;
		x \in \left\{(x_1,\ldots,x_n) \in \R^n : {\textstyle{\sum_{i=1}^{n}}}x_i = 0\right\}.
	\end{equation}
\end{assumption}

\section{Lower Complexity Bounds}\label{sec:lb}

\subsection{Decentralized Subgradient Optimization Algorithms}

In this section, we present the lower bounds on the number of decentralized communications and the number of local subgradient computations required to solve problem~\eqref{eq:main}. These lower bounds apply to a particular class of algorithms, which we refer to as the class of {\em decentralized subgradient optimization algorithms}. This class can be seen as an adaptation of {\em black-box optimization procedures} \citep{scaman2018optimal} to the time-varying network setting, or an adaptation of {\em first-order decentralized optimization algorithms} \citep{kovalev2021lower} to the non-smooth optimization setting.

Non-smooth optimization algorithms typically perform incremental updates by computing the subgradient of a given objective function. The set of all subgradients of a convex function, called the subdifferential, can be multivalued in general. Thus, it is necessary to select the specific subgradient that the algorithm will use. This is done by the {\em subgradient oracle}, which is described by \Cref{def:sub}.
\begin{definition}\label{def:sub}
	For each $i \in \cV$, a function $\subgrad f_i(x)\colon \sX \to \sX$ is called a subgradient oracle associated with the function $f_i(x)$ if, for all $x \in \sX$, it satisfies $\subgrad f_i(x) \in \partial f_i(x)$. That is, for each $i \in \cV$ and for all $x,x' \in \sX$, the following inequality holds:
	\begin{equation}
		f_i(x') \geq f_i(x) + \<\subgrad f_i(x), x'-x>.
	\end{equation}
\end{definition}
Further, we provide the formal description of the class of decentralized subgradient optimization algorithms in the following \Cref{def:alg}.
\begin{definition}\label{def:alg}
	An algorithm is called a decentralized subgradient optimization algorithm with the subgradient computation time $\tloc > 0$ and decentralized communication time $\tcom > 0$ if it satisfies the following constraints:
	\begin{enumerate}
		\item [\bf (i)] {\bf Internal memory.} At any time $\texe \geq 0$, each node $i \in \cV$ maintains an internal memory, which is represented by a set-valued function $\mem{i}(\texe)\colon \R_{+} \to 2^{\sX}$. The internal memory can be updated by subgradient computation or decentralized communication, which is formally represented by the following inclusion:
		      \begin{equation}\label{eq:mem}
			      \mem{i}(\texe) \subset \memloc{i}(\texe) \cup \memcom{i}(\texe),
		      \end{equation}
		      where set-valued functions $ \memloc{i}(\texe),\memcom{i}(\texe)\colon \R_+ \to 2^{\sX}$ are defined below.
		\item [\bf (ii)] {\bf Subgradient computation.} At any time $\texe \geq 0$, each node $i \in \cV$ can update its internal memory $\mem{i}(\texe)$ by computing the subgradient $\subgrad f_i(x)$ of the function $f_i(x)$, which takes time $\tloc$. That is, for all $\texe \geq 0$, the set $\memloc{i}(\texe)$ is defined as follows:
		      \begin{equation}\label{eq:memloc}
			      \memloc{i}(\texe) = \begin{cases}
				      \spanset(\{ x, \subgrad f_i(x) : x \in \mem{i}(\texe - \tloc) \}) & \texe \geq \tloc \\
				      \varnothing                                                       & \texe < \tloc    
			      \end{cases}.
		      \end{equation}
		\item [\bf (iii)] {\bf Decentralized communication.} At any time $\texe \geq 0$, each node $i \in \cV$ can update its internal memory $\mem{i}(\texe)$ by performing decentralized communication across the communication network, which takes time $\tcom$. That is, for all $\texe \geq 0$, the set $\memcom{i}(\texe)$ is defined as follows:
		      \begin{equation}\label{eq:memcom}
			      \memcom{i}(\texe) = \begin{cases}
				      \spanset \bigl(\textstyle{\bigcup}_{(j,i) \in \cE(\texe)} \mem{j}(\texe - \tcom)\bigr) & \texe \geq \tcom \\
				      \varnothing                                                                            & \texe < \tcom
			      \end{cases}.
		      \end{equation}
		\item [\bf (iv)] {\bf Initialization and output.} At time $\texe = 0$, each node $i \in \cV$ must initialize its internal memory with the zero vector, that is, $\mem{i}(0) = \{0\}$. At any time $\texe\geq 0$, each node $i \in \cV$ must specify a single output vector from its internal memory, $x_{o,i}(\texe) \in \mem{i}(\texe)$.
	\end{enumerate}
\end{definition}

\subsection{Lower Bounds}

Now, we are ready to present the lower bounds on the execution time $\texe \geq 0$ required to find an $\epsilon$-approximate solution\footnote{A vector $x \in \sX$ is called an $\epsilon$-approximate solution to problem~\eqref{eq:main} if $p(x) - p(x^*) \leq \epsilon$.} to problem~\eqref{eq:main} by any algorithm satisfying \Cref{def:alg}.
\Cref{thm:lb_sc} provides the lower bound in the strongly convex case ($r>0$), and \Cref{thm:lb} provides the lower bound in the convex case ($r=0$). These lower bounds naturally depend on the precision $\epsilon > 0$, the parameters of the problem, including the Lipschitz constant $M > 0$, the regularization parameter $r \geq 0$, the distance $R > 0$, and the parameters of the network, including the condition number $\chi \geq 1$, communication time $\tcom > 0$, and subgradient computation time $\tloc > 0$.

\begin{theorem}\label{thm:lb_sc}
	For arbitrary parameters $M, r, \epsilon, \tcom, \tloc > 0$ and $\chi \geq 1$, there exists an optimization problem of the form~\eqref{eq:main} satisfying \Cref{ass:cvx,,ass:lip,,ass:sol}, corresponding subgradient oracles given by \Cref{def:sub}, a time varying network $\cG(\texe) = (\cV,\cE(\texe))$, and a corresponding time-varying gossip matrix $\mW(\texe)$ satisfying \Cref{ass:gossip,ass:chi}, such that at least the following time $\texe$ is required to reach precision $p(x_{o,i}(\texe)) - p(x^*) \leq \epsilon$ by any decentralized subgradient optimization algorithm satisfying \Cref{def:alg}:
	\begin{equation}
		\texe \geq \Omega\left(\tcom \cdot \frac{\chi M}{\sqrt{r\epsilon}} + \tloc \cdot \frac{M^2}{r\epsilon}\right).
	\end{equation}
\end{theorem}

\begin{theorem}\label{thm:lb}
	For arbitrary parameters $M, R, \epsilon, \tcom, \tloc > 0$ and $\chi \geq 1$, there exists an optimization problem of the form~\eqref{eq:main} with zero regularization ($r = 0$) satisfying \Cref{ass:cvx,,ass:lip,,ass:sol}, corresponding subgradient oracles given by \Cref{def:sub}, a time varying network $\cG(\texe) = (\cV,\cE(\texe))$, and a corresponding time-varying gossip matrix $\mW(\texe)$ satisfying \Cref{ass:gossip,ass:chi}, such that at least the following time $\texe$ is required to reach precision $p(x_{o,i}(\texe)) - p(x^*) \leq \epsilon$ by any decentralized subgradient optimization algorithm satisfying \Cref{def:alg}:
	\begin{equation}
		\texe \geq \Omega\left(\tcom \cdot \frac{\chi MR}{\epsilon} + \tloc \cdot \frac{M^2R^2}{\epsilon^2}\right).
	\end{equation}
\end{theorem}

The proofs of \Cref{thm:lb_sc,thm:lb} can be found in \Cref{sec:lb_proof}.
Further, we provide a brief and informal description of the main theoretical ideas that underlie these proofs:
\begin{enumerate}
	\item[\bf (i)] We select a specific ``hard'' instance of problem~\eqref{eq:main}. In particular, we choose the objective function of the form $p(x) = a\sum_{j=1}^{d-1} \abs{\<\basis{d}{j+1} - \basis{d}{j},x>} - a\<\basis{d}{1},x>+ \frac{r}{2}\sqn{x}$, which was used by \citet{arjevani2015communication,scaman2018optimal} in the proof of lower bounds on the communication complexity in centralized and fixed-network settings. One can show that the gap $p(x) - p(x^*)$ is lower-bounded by a positive constant as long as the last component of the vector $x$ is zero, and it takes $\Omega(\tloc \cdot d)$ time to break this bound due to the constraint on the subgradient updates~\eqref{eq:memloc}.
	\item[\bf (ii)] We split the objective function between two nodes of a star-topology network with a time-varying central node, which was previously utilized by \citet{kovalev2021lower} in the proof of lower bounds for optimizing smooth functions. One can show that it takes $\Omega(n) = \Omega(\chi)$ communications to exchange information between the two selected nodes due to the time-varying center. This contrasts with the fixed path-topology network used by \citet{scaman2017optimal,scaman2018optimal}, where such an exchange would take $\Omega(n) = \Omega(\sqrt{\chi})$ communications. Moreover, using the constraint~\eqref{eq:memloc}, we can show that it takes $\Omega(\tcom \cdot nd)$ time to make the last component of the vector $x$ nonzero and break the lower bound on the gap $p(x) - p(x^*)$, thanks to the way we split the objective function.
	\item[\bf (iii)] Based on the above considerations, we show that the total execution time required to solve the problem is lower-bounded by $\Omega\left(\tcom \cdot nd + \tloc \cdot d\right)$. Thus, we obtain the desired results by making a specific choice of the dimension $d$, network size $n$, and other parameters of problem~\eqref{eq:main}.
\end{enumerate}

\subsection{Comparison with the Lower Bounds in Centralized and Fixed Network Settings}

\begin{table}[t!]
	\centering
	\caption{Lower bounds on the communication complexity of solving problem~\eqref{eq:main} in the centralized \citep{arjevani2015communication}, decentralized fixed network \citep{scaman2018optimal}, and decentralized time-varying network (\Cref{thm:lb_sc,thm:lb}) settings.}
	\label{tab:lb}
	\begin{NiceTabular}{|c|l|l|l|}
		\CodeBefore
		\Body
		\toprule
		\bf Setting
		 &
		\Block[c]{1-1}{\bf Centralized}
		 &
		\Block[c]{1-1}{\bf Fixed networks\tablefootnote{\citet{scaman2018optimal} do not provide any lower complexity bounds in the strongly convex setting. However, the desired lower bound on the communication complexity can be obtained by extending their analysis.}}
		 &
		\Block[c]{1-1}{\bf Time-varying networks}
		\\\midrule
		\bf Strongly convex
		 & $\Omega\left( M / \sqrt{r\epsilon}\right)$
		 & $\Omega\left( \alertB{\sqrt{\chi}}  M / \sqrt{r\epsilon}\right)$
		 & \hspace{2em}$\Omega\left( \alertB{\chi}  M / \sqrt{r\epsilon}\right)$
		\\\midrule
		\bf Convex
		 & $\Omega\left( MR / \epsilon \right)$
		 & $\Omega\left(  \alertB{\sqrt{\chi}}  MR / \epsilon\right)$
		 & \hspace{2em}$\Omega\left( \alertB{\chi}  MR / \epsilon\right)$
		\\\bottomrule
	\end{NiceTabular}
\end{table}

We compare the lower complexity bounds for solving non-smooth convex optimization problems in the three main distributed optimization settings: centralized, decentralized fixed network, and decentralized time-varying network. The lower subgradient computation complexity bounds coincide in these cases (\citet{nesterov2013introductory},\citet{scaman2018optimal},\Cref{thm:lb_sc,thm:lb}). However, the situation with the communication complexity is different. See \Cref{tab:lb} for a summary.

\Cref{thm:lb_sc,thm:lb} imply that the communication complexity in the decentralized time-varying network setting is proportional to the network condition number $\alertB{\chi}$. In contrast, the communication complexity in the fixed network setting is proportional to $\alertB{\sqrt{\chi}}$, which reflects the fact that time-varying networks are more difficult to deal with compared to fixed networks. In particular, there was a long-standing conjecture that the ``upgrade'' from the factor $\alertB{\chi}$ to the factor $\alertB{\sqrt{\chi}}$ in communication complexity is impossible in the time-varying network setting. Only recently, this conjecture was proved for smooth functions by \citet{kovalev2021lower}, and now we resolve this open question in the non-smooth case as well.

\section{Optimal Algorithm}\label{sec:alg}

In this section, we develop an optimal algorithm for solving the non-smooth convex decentralized optimization problem~\eqref{eq:main} over time-varying networks. The design of our algorithm relies on a specific saddle-point reformulation of the problem, which we describe in the following section.

\subsection{Saddle-Point Reformulation}

Let functions $F(x) \colon (\sX)^n \to \R$ and $G(y,z)\colon (\sX)^n\times (\sX)^n \to \R$ be defined as follows:
\begin{equation}\label{eq:FG}
	F(x) = \sum_{i=1}^{n} f_i(x_i) + \frac{r_x}{2}\sqn{x}
	\quad
	\text{and}
	\quad
	G(y,z) = \frac{r_{yz}}{2}\sqn{y+z},
\end{equation}
where $x = (x_1,\ldots,x_n) \in (\sX)^n$, and $r_x,r_{yz} > 0$ are some constants that satisfy
\begin{equation}\label{eq:r_xyz}
	r_x + 1/r_{yz} = r.
\end{equation}
Consider the following saddle-point problem:
\begin{equation}\label{eq:spp}
	\adjustlimits\min_{x \in (\sX)^n}\max_{y \in (\sX)^n}\max_{z \in (\sX)^n} \left[Q(x,y,z) = F(x) - \<y,x> - G(y,z)\right]
	\quad
	\text{s.t.}
	\quad
	z \in \cL^\perp,
\end{equation}
where $\cL^\perp \subset (\sX)^n$ is the orthogonal complement to the so-called consensus space $\cL \subset (\sX)^n$, defined as follows:
\begin{equation}\label{eq:consensus_space}
	\cL = \{(x_1,\ldots,x_n) : x_1=\ldots = x_n\},
	\quad
	\cL^\perp =  \{(x_1,\ldots,x_n) : {\textstyle\sum_{i=1}^n} x_i= 0\}.
\end{equation}
One can show that the saddle-point problem~\eqref{eq:spp} is equivalent to the minimization problem~\eqref{eq:main}. This is justified by the following \Cref{lem:spp}. The proof of the lemma can be found in the \Cref{sec:proof:spp}.
\begin{lemma}\label{lem:spp}
	Problem~\eqref{eq:spp} is equivalent to problem~\eqref{eq:main} in the following sense:
	\begin{equation}
		\adjustlimits\min_{x \in (\sX)^n}\max_{y \in (\sX)^n}\max_{z \in \cL^\perp}Q(x,y,z) = n\cdot \min_{x \in \sX} p(x).
	\end{equation}
\end{lemma}
The saddle-point reformulation of the form~\eqref{eq:spp} was first introduced by \citet{kovalev2020optimal,kovalev2021lower} to develop optimal decentralized algorithms for optimizing smooth functions. However, these are not applicable to the non-smooth case. To the best of our knowledge, the only attempt to adapt the reformulation~\eqref{eq:spp} to the non-smooth setting was made by \citet{lobanov2023non}. However, their results have significant downsides, which we discuss in \Cref{sec:alg_compare}.

\subsection{New Algorithm and its Convergence}

\begin{algorithm}[t]
	\caption{ }
	\label{alg}
	\begin{algorithmic}[1]
		\State {\bf input:} $x^0\in (\sX)^n$, $y^0\in (\sX)^n$, $z^0\in \cL^\perp$, $m^0 \in (\sX)^n$
		\State {\bf parameters:}
		$K,T \in \{1,2,\ldots\}$, $\{(\alpha_k,\beta_k,\gamma_k,\sigma_k,\lambda_k,\tau_x^k,\eta_x^k,\eta_y^k,\eta_z^k, \theta_z^k)\}_{k=0}^{K-1}\subset \R_+^{10}$
		\For{$k = 0,1,\ldots,K-1$}
		\State $\uln{y}^k = \alpha_k y^k + (1-\alpha_k)\ol{y}^k$,\quad
		$\uln{z}^k = \alpha_k z^k + (1-\alpha_k)\ol{z}^k$
		\label{line:comb}
		\State $g_y^k = \nabla_y G(\uln{y}^k,\uln{z}^k)$,\quad
		$g_z^k = \nabla_z G(\uln{y}^k,\uln{z}^k)$,\quad
		where function $G(y,z)$ is defined in \cref{eq:FG}
		\label{line:grad}
		\State $\tilde{g}_z^k = (\mW_k \otimes\mI_d)g_z^k$,\;\;
		$\hat{g}_z^k =  (\mW_k \otimes\mI_d)(g_z^k + m^k)$,\\
		\hspace{\algorithmicindent}where $\mW_k$ denotes the gossip matrix $\mW(\texe)$ at the current time $\texe$
		\label{line:grad_comm}
		\State $y^{k+1} = y^k - \eta_y^k (g_y^k + \hat{x}^{k+1})$,\quad
		$z^{k+1} = z^k - \eta_z^k \hat{g}_z^k$,\quad
		$\hat{x}^{k+1} = x^k + \gamma_k (\tilde{x}^k - x^{k-1})$
		\label{line:yz}
		\State $\ol{y}^{k+1} = \uln{y}^k + \alpha_k(y^{k+1} - y^k)$,\quad
		$\ol{z}^{k+1} = \uln{z}^k - \theta_z^k \tilde{g}_z^k$,\quad
		$m^{k+1} = (\eta_z^k/\eta_z^{k+1})(m^k +g_z^k  - \hat{g}_z^k)$
		\label{line:ext}
		\State $x^{k,0} = x^k$ \label{line:iinit}
		\For{$t = 0,1,\ldots,T-1$}
		\State $g_x^{k,t} = (\subgrad f_1(x_1^{k,t}),\ldots,\subgrad f_n(x_n^{k,t}))$ \label{line:ig}
		\State $x^{k,t+1} = x^{k,t} - \eta_x^k \left(g_x^{k,t} + \beta_k x^{k,t+1} - y^{k+1} + \tau_x^k(x^{k,t+1} - x^k) \right)$ \label{line:ix}
		\EndFor
		\State $x^{k+1} = \sigma_k x^{k,T} + (1-\sigma_k)\tilde{x}^{k+1}$,\quad
		$\tilde{x}^{k+1} = \frac{1}{T}\sum_{t=1}^{T} x^{k,t}$,\quad
		$\ol{x}^{k+1} = \alpha_k \tilde{x}^{k+1} + (1-\alpha_k)\ol{x}^k$
		\label{line:ox}
		\EndFor
		\State $(x_a^K,y_a^K,z_a^K) = (\sum_{k=1}^K \lambda_k)^{-1}\sum_{k=1}^K \lambda_k(\ol{x}^k,\ol{y}^k,\ol{z}^k)$
		\label{line:avg}
		\State {\bf output:} $x_o^K = \frac{1}{n}\sum_{i=1}^{n}x_{a,i}^K \in \sX$,\quad
		where\;$(x_{a,1}^K,\ldots,x_{a,n}^K) = x_a^K \in (\sX)^n$
		\label{line:output}
	\end{algorithmic}
\end{algorithm}

Now, we present \Cref{alg} for solving problem~\eqref{eq:main}. We provide upper bounds on the number of decentralized communications $K$ and the number of subgradient computations $K\times T$ required to find an $\epsilon$-approximate solution to the problem. \Cref{thm:alg_sc,,thm:alg} provide the upper bounds in the strongly convex ($r > 0$) and convex ($r = 0$) cases, respectively. The proofs can be found in \Cref{sec:alg_proof}. The total execution time of \Cref{alg} is upper-bounded as $\texe = \cO \left(\tcom \cdot K + \tloc \cdot K\times T\right)$, where the communication time $\tcom > 0$ and the subgradient computation time $\tloc > 0$ are described in \Cref{def:alg}. This upper-bound on the execution time cannot be improved because of the lower bounds established in the previous \Cref{sec:lb}. Therefore, \Cref{alg} is an optimal algorithm for solving problem~\eqref{eq:main}.

\begin{theorem}\label{thm:alg_sc}
	Under \Cref{ass:cvx,,ass:lip,,ass:sol,,ass:gossip,,ass:chi}, let $r > 0$ (strongly convex case). Then \Cref{alg} requires $K = \cO\left( \frac{\chi M}{\sqrt{r\epsilon}}\right)$ decentralized communications (\cref{line:grad_comm} of \Cref{alg}) and $ K\times T = \cO\left(\frac{M^2}{r\epsilon}\right)$ subgradient computations (\cref{line:ig} of \Cref{alg}) to reach precision $p(x_o^K) - p(x^*) \leq \epsilon$.
\end{theorem}

\begin{theorem}\label{thm:alg}
	Under \Cref{ass:cvx,,ass:lip,,ass:sol,,ass:gossip,,ass:chi}, let $r = 0$ (convex case). Then \Cref{alg} requires $K = \cO\left( \frac{\chi MR}{\epsilon}\right)$ decentralized communications (\cref{line:grad_comm} of \Cref{alg}) and $K\times T = \cO\left(\frac{M^2 R^2}{\epsilon^2}\right)$ subgradient computations (\cref{line:ig} of \Cref{alg}) to reach precision $p(x_o^K) - p(x^*) \leq \epsilon$.
\end{theorem}

The design of \Cref{alg} is based on the fundamental Forward-Backward algorithm \citep{bauschke2011convex}. Let $\mathsf{E} = (\sX)^n \times (\sX)^n \times \cL^\perp$ be a Euclidean space, and consider a monotone operator $A(u)\colon \mathsf{E} \to \mathsf{E}$ and a maximally-monotone multivalued operator $B(u)\colon \mathsf{E} \to 2^\mathsf{E}$ defined as follows:
\begin{equation}
	A(u) =
	\left[\begin{NiceMatrix}[r]
			\Block[c]{1-1}{0} \\\nabla_y G(y,z)\\\mP \nabla_z G(y,z)
		\end{NiceMatrix}\right],
	\qquad
	B(u) =
	\left[\begin{NiceMatrix}[c]
			\partial F(x) - y \\ x\\0
		\end{NiceMatrix}\right],
\end{equation}
where $u = (x,y,z) \in \mathsf{E}$, and $\mP = (\mI_n - (1/n)\ones_n\ones_n^\top)\otimes \mI_d \in \R^{nd\times nd}$ is the orthogonal projection matrix onto $\cL^\perp$. Then problem~\eqref{eq:spp} is equivalent to the following monotone inclusion problem:
\begin{equation}\label{eq:inclusion}
	\text{find } u \in \mathsf{E} \text{ such that } 0 \in A(u) + B(u).
\end{equation}
The basic Forward-Backward algorithm iterates $u^{k+1} = (\mathrm{id} + B)^{-1}(u^k - A(u^k))$, where $\mathrm{id}$ is the identity operator and $(\mathrm{id} + B)^{-1}$ denotes the inverse of the operator $\mathrm{id}(u) + B(u)$, which is called resolvent. \Cref{alg} can be obtained by making the following major modifications to these iterations:
\begin{enumerate}
	\item[\bf (i)] We accelerate the convergence of the Forward-Backward algorithm using Nesterov acceleration \citep{nesterov1983method}. Although this mechanism cannot be applied to the general monotone inclusion problem~\eqref{eq:inclusion}, \citet{kovalev2020optimal} showed that it can be used when the operator $A(u)$ is equal to the gradient of a smooth convex function, which is true in our case.
	\item[\bf (ii)] Computation of the operator $A(u)$ requires multiplication with the matrix $\mP$. This, in turn, requires an exact averaging of a vector, which is difficult to do over the time-varying network. \citet{kovalev2021adom} showed that this obstacle can be tackled with the Error-Feedback mechanism for decentralized communication, which we also utilize.
	\item[\bf (iii)] At each iteration of the algorithm, we have to compute the resolvent, which requires solving an auxiliary subproblem $\min_{x}\max_{y} \frac{\tau_x}{2}\sqn{x-x^k} + F(x) -\<y,x> - \frac{\tau_y}{2}\sqn{y-y^k}$. This problem cannot be solved exactly, so we have to find an approximate solution using an additional ``inner'' algorithm based on the subgradient method \citep{nesterov2013introductory} and the Chambolle-Pock operator splitting \citep{chambolle2011first}. We also have to conduct a careful analysis to find an efficient way to combine the inner and the ``outer'' Forward-Backward algorithms and avoid unnecessary waste of subgradient calls.
\end{enumerate}
The design of \Cref{alg} shares some similarities with the algorithm of \citet{kovalev2021lower} such as {\bf(i)} and {\bf(ii)} above. However, \citet{kovalev2021lower} simply add the gradient $\nabla F(x)$ to the operator $A(u)$ and use the accelerated version of the Forward-Backward algorithm, which we obviously cannot do as the function $F(x)$ is not smooth. Instead, we have to put the subdifferential $\partial F(x)$ into the operator $B(u)$ and follow {\bf (iii)} above. Part {\bf(iii)}, in turn, shares some similarities with the algorithm of \citet{lan2020communication}. However, \citet{lan2020communication} simply have a zero operator $A(u) = 0$, which makes {\bf (i)} and {\bf (ii)} above unnecessary in their case. In contrast, we cannot make such simplifications because we work in the much more complicated setting of time-varying networks.

\subsection{Comparison with the Existing Results}\label{sec:alg_compare}

\begin{table}[t]
	\centering
	\caption{The execution time $\texe$ required to find an $\epsilon$-approximate solution to the decentralized optimization problem~\eqref{eq:main} by the following algorithms: \algname{D-SubGD} \citep{nedic2009distributed}, \algname{SubGD-Push} \citep{nedic2014distributed}, \algname{ZO-SADOM} \citep{lobanov2023non}, and \Cref{alg} (this paper). \sethlcolor{green!40}\hl{Decentralized communication} and \sethlcolor{yellow!90}\hl{subgradient computation} complexities are marked with \sethlcolor{green!40}\hl{green} and \sethlcolor{yellow!90}\hl{yellow} colors, respectively. For \algname{D-SubGD}, the complexity is not provided because the algorithm converges only to a neighborhood of the solution. For \algname{SubGD-Push}, $\mathrm{poly}(M,R,d)$ denotes a certain polynomial in $M,R,d$. For \algname{ZO-SADOM}, the differences from the optimal complexities are highlighted in \alertR{red color}. }
	\label{tab:alg}
	\resizebox{\textwidth}{!}{
		\begin{NiceTabular}{|c|c|c|}
			\CodeBefore
			\cellcolor{cyan!12}{5-1,5-2,5-3,6-1,6-2,6-3}
			\Body
			\toprule
			\bf Algorithm
			 &
			\bf Strongly-convex case complexity
			 &
			\bf Convex case complexity
			\\\midrule
			\algname{D-SubGD}
			 &
			\Block{1-2}{N/A}
			 &
			\\\midrule
			\algname{SubGD-Push}
			 &
			\Block{1-2}{
				{$\displaystyle
							\tcom\cdot\hlcom{\frac{\mathrm{poly}(M,R,d)\cdot n^{2n}\log^2\tfrac{1}{\epsilon}}{\epsilon^2}}
							+
							\tloc \cdot \hlloc{\frac{\mathrm{poly}(M,R,d)\cdot n^{2n}\log^2\tfrac{1}{\epsilon}}{\epsilon^2}}$}
			}
			 &
			\\\midrule
			\algname{ZO-SADOM}
			 &
			{$\displaystyle
						\tcom\cdot\hlcom{\frac{\chi M\alertR{d^{1/4}\log\tfrac{1}{\epsilon}}}{\sqrt{r\epsilon}}}
						+
						\tloc \cdot \hlloc{\frac{M^2\alertR{d\log\tfrac{1}{\epsilon}}}{r\epsilon}}$}
			 &
			{$\displaystyle
						\tcom\cdot\hlcom{\frac{\chi MR\alertR{d^{1/4}\log\tfrac{1}{\epsilon}}}{\epsilon}}
						+
						\tloc \cdot \hlloc{\frac{M^2R^2\alertR{d\log\tfrac{1}{\epsilon}}}{\epsilon^2}}$}
			\\\midrule
			\makecell{\bf \Cref{alg}}
			 &
			{$\displaystyle
						\tcom \cdot \hlcom{\frac{\chi M}{\sqrt{r\epsilon}}}
						+
						\tloc \cdot \hlloc{\frac{M^2}{r\epsilon}}$}
			 &
			{$\displaystyle
						\tcom \cdot \hlcom{\frac{\chi MR}{\epsilon}}
						+
						\tloc \cdot \hlloc{\frac{M^2R^2}{\epsilon^2}}$}
			\\\midrule
			\makecell{\bf Lower Bounds}
			 &
			{$\displaystyle
						\tcom \cdot \hlcom{\frac{\chi M}{\sqrt{r\epsilon}}}
						+
						\tloc \cdot \hlloc{\frac{M^2}{r\epsilon}}$}
			 &
			{$\displaystyle
						\tcom \cdot \hlcom{\frac{\chi MR}{\epsilon}}
						+
						\tloc \cdot \hlloc{\frac{M^2R^2}{\epsilon^2}}$}
			\\\bottomrule
		\end{NiceTabular}}
\end{table}

One could naturally expect that the existing optimal algorithms, originally developed for fixed networks, such as \algname{DCS} \citep{lan2020communication} and \algname{MSPD} \citep{scaman2018optimal}, could be applied to solve problem~\eqref{eq:main} over time-varying networks. However, this is not the case, which is justified by the lack of corresponding theoretical guarantees and was shown empirically by \citet{kovalev2021adom}. Therefore, we have to consider only those algorithms that were specifically developed for the time-varying network setting.

We provide a comparison of our \Cref{alg} with the existing state-of-the-art decentralized methods for solving convex non-smooth optimization problems over time-varying networks in \Cref{tab:alg}.\footnote{We ignore universal constants in \Cref{tab:alg} like in the $\cO(\cdot)$ and $\Omega(\cdot)$ notation.} These include \algname{D-SubGD} \citep{nedic2009distributed}, \algname{SubGD-Push} \citep{nedic2014distributed}, and \algname{ZO-SADOM} \citep{lobanov2023non}. The first two algorithms have poor performance: \algname{D-SubGD} converges only to limited precision, and \algname{SubGD-Push} converges at a slow rate of $\cO(\log^2(1/\epsilon)/\epsilon^2)$, which does not match even the iteration complexity of the standard centralized subgradient method, let alone the improved complexity of \Cref{alg}. The complexity of \algname{ZO-SADOM} is also worse than the lower bounds. Moreover, the theoretical results of \citet{lobanov2023non} have substantial drawbacks compared to ours:
\begin{enumerate}
	\item[\bf (i)] \citet{lobanov2023non} do not provide any theoretical insights or innovations in the analysis of their algorithm. In particular, they use the randomized smoothing technique \citep{duchi2012randomized} to obtain a smooth approximation of the objective $p(x)$, and apply the existing algorithm of \citet{kovalev2021lower} to minimize this approximation. In contrast, we develop a new algorithm that directly works with the original non-smooth objective $p(x)$.
	\item[\bf(ii)] \algname{ZO-SADOM} has extra factors $\alertR{d^{1/4}\log(1/\epsilon)}$ and $\alertR{d\log(1/\epsilon)}$ in the decentralized communication and subgradient computation complexities, respectively, compared to the optimal complexity of our \Cref{alg}. Thus, the performance of \algname{ZO-SADOM} can be poor when applied, for instance, to large-scale machine learning problems in which the dimension $d$ can be huge.
\end{enumerate}


\bibliographystyle{apalike}
\bibliography{ref.bib}

\begin{thebibliography}{}

\bibitem[Allen-Zhu, 2018]{allen2018natasha}
Allen-Zhu, Z. (2018).
\newblock Natasha 2: Faster non-convex optimization than sgd.
\newblock {\em Advances in neural information processing systems}, 31.

\bibitem[Arjevani and Shamir, 2015]{arjevani2015communication}
Arjevani, Y. and Shamir, O. (2015).
\newblock Communication complexity of distributed convex learning and optimization.
\newblock {\em Advances in neural information processing systems}, 28.

\bibitem[Bauschke and Combettes, 2011]{bauschke2011convex}
Bauschke, H.~H. and Combettes, P.~L. (2011).
\newblock {\em Convex Analysis and Monotone Operator Theory in Hilbert Spaces}.
\newblock Springer Science \& Business Media.

\bibitem[Bazerque and Giannakis, 2009]{bazerque2009distributed}
Bazerque, J.~A. and Giannakis, G.~B. (2009).
\newblock Distributed spectrum sensing for cognitive radio networks by exploiting sparsity.
\newblock {\em IEEE Transactions on Signal Processing}, 58(3):1847--1862.

\bibitem[Beck et~al., 2014]{beck20141}
Beck, A., Nedi{\'c}, A., Ozdaglar, A., and Teboulle, M. (2014).
\newblock An $ o (1/k) $ gradient method for network resource allocation problems.
\newblock {\em IEEE Transactions on Control of Network Systems}, 1(1):64--73.

\bibitem[Carmon et~al., 2018]{carmon2018accelerated}
Carmon, Y., Duchi, J.~C., Hinder, O., and Sidford, A. (2018).
\newblock Accelerated methods for nonconvex optimization.
\newblock {\em SIAM Journal on Optimization}, 28(2):1751--1772.

\bibitem[Carmon et~al., 2020]{carmon2020lower}
Carmon, Y., Duchi, J.~C., Hinder, O., and Sidford, A. (2020).
\newblock Lower bounds for finding stationary points i.
\newblock {\em Mathematical Programming}, 184(1):71--120.

\bibitem[Chambolle and Pock, 2011]{chambolle2011first}
Chambolle, A. and Pock, T. (2011).
\newblock A first-order primal-dual algorithm for convex problems with applications to imaging.
\newblock {\em Journal of mathematical imaging and vision}, 40:120--145.

\bibitem[Duchi et~al., 2011]{duchi2011adaptive}
Duchi, J., Hazan, E., and Singer, Y. (2011).
\newblock Adaptive subgradient methods for online learning and stochastic optimization.
\newblock {\em Journal of machine learning research}, 12(7).

\bibitem[Duchi et~al., 2012]{duchi2012randomized}
Duchi, J.~C., Bartlett, P.~L., and Wainwright, M.~J. (2012).
\newblock Randomized smoothing for stochastic optimization.
\newblock {\em SIAM Journal on Optimization}, 22(2):674--701.

\bibitem[Gan et~al., 2012]{gan2012optimal}
Gan, L., Topcu, U., and Low, S.~H. (2012).
\newblock Optimal decentralized protocol for electric vehicle charging.
\newblock {\em IEEE Transactions on Power Systems}, 28(2):940--951.

\bibitem[Giselsson et~al., 2013]{giselsson2013accelerated}
Giselsson, P., Doan, M.~D., Keviczky, T., De~Schutter, B., and Rantzer, A. (2013).
\newblock Accelerated gradient methods and dual decomposition in distributed model predictive control.
\newblock {\em Automatica}, 49(3):829--833.

\bibitem[Hendrikx et~al., 2021]{hendrikx2021optimal}
Hendrikx, H., Bach, F., and Massoulie, L. (2021).
\newblock An optimal algorithm for decentralized finite-sum optimization.
\newblock {\em SIAM Journal on Optimization}, 31(4):2753--2783.

\bibitem[Hinton et~al., 2012]{hinton2012neural}
Hinton, G., Srivastava, N., and Swersky, K. (2012).
\newblock Neural networks for machine learning lecture 6a overview of mini-batch gradient descent.
\newblock {\em Cited on}, 14(8):2.

\bibitem[Karimi et~al., 2016]{karimi2016linear}
Karimi, H., Nutini, J., and Schmidt, M. (2016).
\newblock Linear convergence of gradient and proximal-gradient methods under the polyak-{\l}ojasiewicz condition.
\newblock In {\em Machine Learning and Knowledge Discovery in Databases: European Conference, ECML PKDD 2016, Riva del Garda, Italy, September 19-23, 2016, Proceedings, Part I 16}, pages 795--811. Springer.

\bibitem[Karimireddy et~al., 2020]{karimireddy2020scaffold}
Karimireddy, S.~P., Kale, S., Mohri, M., Reddi, S., Stich, S., and Suresh, A.~T. (2020).
\newblock Scaffold: Stochastic controlled averaging for federated learning.
\newblock In {\em International conference on machine learning}, pages 5132--5143. PMLR.

\bibitem[Kingma and Ba, 2014]{kingma2014adam}
Kingma, D.~P. and Ba, J. (2014).
\newblock Adam: A method for stochastic optimization.
\newblock {\em arXiv preprint arXiv:1412.6980}.

\bibitem[Kolar et~al., 2010]{kolar2010estimating}
Kolar, M., Song, L., Ahmed, A., and Xing, E.~P. (2010).
\newblock Estimating time-varying networks.
\newblock {\em The Annals of Applied Statistics}, pages 94--123.

\bibitem[Konecn{\`y} et~al., 2016]{konecny2016federated}
Konecn{\`y}, J., McMahan, H.~B., Yu, F.~X., Richt{\'a}rik, P., Suresh, A.~T., and Bacon, D. (2016).
\newblock Federated learning: Strategies for improving communication efficiency.
\newblock {\em arXiv preprint arXiv:1610.05492}, 8.

\bibitem[Kovalev et~al., 2022]{kovalev2022optimal}
Kovalev, D., Beznosikov, A., Sadiev, A., Persiianov, M., Richt{\'a}rik, P., and Gasnikov, A. (2022).
\newblock Optimal algorithms for decentralized stochastic variational inequalities.
\newblock {\em Advances in Neural Information Processing Systems}, 35:31073--31088.

\bibitem[Kovalev et~al., 2021a]{kovalev2021lower}
Kovalev, D., Gasanov, E., Gasnikov, A., and Richtarik, P. (2021a).
\newblock Lower bounds and optimal algorithms for smooth and strongly convex decentralized optimization over time-varying networks.
\newblock {\em Advances in Neural Information Processing Systems}, 34:22325--22335.

\bibitem[Kovalev et~al., 2020]{kovalev2020optimal}
Kovalev, D., Salim, A., and Richt{\'a}rik, P. (2020).
\newblock Optimal and practical algorithms for smooth and strongly convex decentralized optimization.
\newblock {\em Advances in Neural Information Processing Systems}, 33:18342--18352.

\bibitem[Kovalev et~al., 2021b]{kovalev2021adom}
Kovalev, D., Shulgin, E., Richt{\'a}rik, P., Rogozin, A.~V., and Gasnikov, A. (2021b).
\newblock Adom: Accelerated decentralized optimization method for time-varying networks.
\newblock In {\em International Conference on Machine Learning}, pages 5784--5793. PMLR.

\bibitem[Lan et~al., 2020]{lan2020communication}
Lan, G., Lee, S., and Zhou, Y. (2020).
\newblock Communication-efficient algorithms for decentralized and stochastic optimization.
\newblock {\em Mathematical Programming}, 180(1):237--284.

\bibitem[Li et~al., 2018]{li2018sharp}
Li, H., Fang, C., Yin, W., and Lin, Z. (2018).
\newblock A sharp convergence rate analysis for distributed accelerated gradient methods.
\newblock {\em arXiv preprint arXiv:1810.01053}.

\bibitem[Li and Lin, 2021]{li2021accelerated}
Li, H. and Lin, Z. (2021).
\newblock Accelerated gradient tracking over time-varying graphs for decentralized optimization.
\newblock {\em arXiv preprint arXiv:2104.02596}.

\bibitem[Li et~al., 2022]{li2022variance}
Li, H., Lin, Z., and Fang, Y. (2022).
\newblock Variance reduced extra and diging and their optimal acceleration for strongly convex decentralized optimization.
\newblock {\em Journal of Machine Learning Research}, 23(222):1--41.

\bibitem[Li et~al., 2019]{li2019decentralized}
Li, Z., Shi, W., and Yan, M. (2019).
\newblock A decentralized proximal-gradient method with network independent step-sizes and separated convergence rates.
\newblock {\em IEEE Transactions on Signal Processing}, 67(17):4494--4506.

\bibitem[Lian et~al., 2017]{lian2017can}
Lian, X., Zhang, C., Zhang, H., Hsieh, C.-J., Zhang, W., and Liu, J. (2017).
\newblock Can decentralized algorithms outperform centralized algorithms? a case study for decentralized parallel stochastic gradient descent.
\newblock {\em Advances in neural information processing systems}, 30.

\bibitem[Lobanov et~al., 2023]{lobanov2023non}
Lobanov, A., Veprikov, A., Konin, G., Beznosikov, A., Gasnikov, A., and Kovalev, D. (2023).
\newblock Non-smooth setting of stochastic decentralized convex optimization problem over time-varying graphs.
\newblock {\em Computational Management Science}, 20(1):48.

\bibitem[McMahan et~al., 2017]{mcmahan2017communication}
McMahan, B., Moore, E., Ramage, D., Hampson, S., and y~Arcas, B.~A. (2017).
\newblock Communication-efficient learning of deep networks from decentralized data.
\newblock In {\em Artificial intelligence and statistics}, pages 1273--1282. PMLR.

\bibitem[Metelev et~al., 2024]{metelev2024decentralized}
Metelev, D., Chezhegov, S., Rogozin, A., Kovalev, D., Beznosikov, A., Sholokhov, A., and Gasnikov, A. (2024).
\newblock Decentralized finite-sum optimization over time-varying networks.
\newblock {\em arXiv preprint arXiv:2402.02490}.

\bibitem[Mishchenko et~al., 2022]{mishchenko2022proxskip}
Mishchenko, K., Malinovsky, G., Stich, S., and Richt{\'a}rik, P. (2022).
\newblock Proxskip: Yes! local gradient steps provably lead to communication acceleration! finally!
\newblock In {\em International Conference on Machine Learning}, pages 15750--15769. PMLR.

\bibitem[Necoara et~al., 2019]{necoara2019linear}
Necoara, I., Nesterov, Y., and Glineur, F. (2019).
\newblock Linear convergence of first order methods for non-strongly convex optimization.
\newblock {\em Mathematical Programming}, 175:69--107.

\bibitem[Nedi{\'c} and Olshevsky, 2014]{nedic2014distributed}
Nedi{\'c}, A. and Olshevsky, A. (2014).
\newblock Distributed optimization over time-varying directed graphs.
\newblock {\em IEEE Transactions on Automatic Control}, 60(3):601--615.

\bibitem[Nedic et~al., 2017]{nedic2017achieving}
Nedic, A., Olshevsky, A., and Shi, W. (2017).
\newblock Achieving geometric convergence for distributed optimization over time-varying graphs.
\newblock {\em SIAM Journal on Optimization}, 27(4):2597--2633.

\bibitem[Nedic and Ozdaglar, 2009]{nedic2009distributed}
Nedic, A. and Ozdaglar, A. (2009).
\newblock Distributed subgradient methods for multi-agent optimization.
\newblock {\em IEEE Transactions on Automatic Control}, 54(1):48--61.

\bibitem[Nesterov, 1983]{nesterov1983method}
Nesterov, Y. (1983).
\newblock A method for unconstrained convex minimization problem with the rate of convergence o (1/k2).
\newblock In {\em Dokl. Akad. Nauk. SSSR}, volume 269, page 543.

\bibitem[Nesterov, 2009]{nesterov2009primal}
Nesterov, Y. (2009).
\newblock Primal-dual subgradient methods for convex problems.
\newblock {\em Mathematical programming}, 120(1):221--259.

\bibitem[Nesterov, 2013]{nesterov2013introductory}
Nesterov, Y. (2013).
\newblock {\em Introductory lectures on convex optimization: A basic course}, volume~87.
\newblock Springer Science \& Business Media.

\bibitem[Rabbat and Nowak, 2004]{rabbat2004distributed}
Rabbat, M. and Nowak, R. (2004).
\newblock Distributed optimization in sensor networks.
\newblock In {\em Proceedings of the 3rd international symposium on Information processing in sensor networks}, pages 20--27.

\bibitem[Rogozin et~al., 2021]{rogozin2021towards}
Rogozin, A., Lukoshkin, V., Gasnikov, A., Kovalev, D., and Shulgin, E. (2021).
\newblock Towards accelerated rates for distributed optimization over time-varying networks.
\newblock In {\em Optimization and Applications: 12th International Conference, OPTIMA 2021, Petrovac, Montenegro, September 27--October 1, 2021, Proceedings 12}, pages 258--272. Springer.

\bibitem[Ryabinin et~al., 2021]{ryabinin2021moshpit}
Ryabinin, M., Gorbunov, E., Plokhotnyuk, V., and Pekhimenko, G. (2021).
\newblock Moshpit sgd: Communication-efficient decentralized training on heterogeneous unreliable devices.
\newblock {\em Advances in Neural Information Processing Systems}, 34:18195--18211.

\bibitem[Ryabinin and Gusev, 2020]{ryabinin2020towards}
Ryabinin, M. and Gusev, A. (2020).
\newblock Towards crowdsourced training of large neural networks using decentralized mixture-of-experts.
\newblock {\em Advances in Neural Information Processing Systems}, 33:3659--3672.

\bibitem[Sadiev et~al., 2022]{sadiev2022communication}
Sadiev, A., Kovalev, D., and Richt{\'a}rik, P. (2022).
\newblock Communication acceleration of local gradient methods via an accelerated primal-dual algorithm with an inexact prox.
\newblock {\em Advances in Neural Information Processing Systems}, 35:21777--21791.

\bibitem[Scaman et~al., 2017]{scaman2017optimal}
Scaman, K., Bach, F., Bubeck, S., Lee, Y.~T., and Massouli{\'e}, L. (2017).
\newblock Optimal algorithms for smooth and strongly convex distributed optimization in networks.
\newblock In {\em international conference on machine learning}, pages 3027--3036. PMLR.

\bibitem[Scaman et~al., 2018]{scaman2018optimal}
Scaman, K., Bach, F., Bubeck, S., Massouli{\'e}, L., and Lee, Y.~T. (2018).
\newblock Optimal algorithms for non-smooth distributed optimization in networks.
\newblock {\em Advances in Neural Information Processing Systems}, 31.

\bibitem[Scutari and Sun, 2019]{scutari2019distributed}
Scutari, G. and Sun, Y. (2019).
\newblock Distributed nonconvex constrained optimization over time-varying digraphs.
\newblock {\em Mathematical Programming}, 176:497--544.

\bibitem[Shi et~al., 2015]{shi2015extra}
Shi, W., Ling, Q., Wu, G., and Yin, W. (2015).
\newblock Extra: An exact first-order algorithm for decentralized consensus optimization.
\newblock {\em SIAM Journal on Optimization}, 25(2):944--966.

\bibitem[Zadeh, 1961]{zadeh1961time}
Zadeh, L.~A. (1961).
\newblock Time-varying networks, i.
\newblock {\em Proceedings of the IRE}, 49(10):1488--1503.

\end{thebibliography}


\appendix
\newpage


\part*{Appendix}

\section{Proof of Lemma~\ref{lem:spp}}\label{sec:proof:spp}

The orthogonal complement $\cL^\perp$ to the consensus space $\cL$ is given as follows:
\begin{equation}\label{eq:dual_space}
	\cL^\perp = \left\{(x_1,\ldots,x_n) \in (\sX)^n : x_1+\ldots + x_n = 0\right\}.
\end{equation}
Let us perform the maximization of $Q(x,y,z)$ in the variable $y \in (\sX)^n$:
\begin{align*}
	\max_{y \in (\sX)^n} Q(x,y,z)
	 & \aeq{uses the definition of $Q(x,y,z)$ \cref{eq:spp}}
	\max_{y \in (\sX)^n}F(x) + \<y,x>-G(y,z)
	\\&\aeq{uses the definition of $G(y,z)$ in \cref{eq:FG}}
	F(x) + \max_{y \in (\sX)^n} \left[\<y,x>-\frac{r_{yz}}{2}\sqn{y+z}\right]
	\\&=
	F(x) + \frac{1}{2r_{yz}}\sqn{x} - \<x,z>,
\end{align*}
where \annotate. Next, we perform maximization in the variable $z \in \cL^\perp$:
\begin{align*}
	\adjustlimits\max_{z \in \cL^\perp}\max_{y \in (\sX)^n} Q(x,y,z)
	 & =
	\max_{z \in \cL^\perp} \left[F(x) + \frac{1}{2r_{yz}}\sqn{x} - \<x,z>\right]
	\\&=
	F(x) + \frac{1}{2r_{yz}}\sqn{x} + \max_{z \in \cL^\perp} \left[ - \<x,z>\right]
	\\&=
	F(x) + \frac{1}{2r_{yz}}\sqn{x} + I_{\cL}(x),
\end{align*}
where $I_{\cL}(x)\colon (\sX)^n \to \R$ is the indicator function, which is defined as follows:
\begin{equation}
	I_{\cL}(x) =\max_{z \in \cL^\perp} \left[ - \<x,z>\right]= \begin{cases}
		0       & x \in \cL        \\
		+\infty & \text{otherwise}
	\end{cases}.
\end{equation}
Now, we can rewrite the saddle-point problem~\eqref{eq:spp} as follows
\begin{align*}
	\adjustlimits\min_{x \in (\sX)^n}\max_{y \in (\sX)^n}\max_{z \in \cL^\perp} Q(x,y,z)
	 & \aeq{uses the fact that we can exhange the order of the two consecutive maximizations}
	\adjustlimits\min_{x \in (\sX)^n}\max_{z \in \cL^\perp}\max_{y \in (\sX)^n} Q(x,y,z)
	\\&\aeq{uses the previous equation}
	\min_{x \in (\sX)^n}F(x) + \frac{1}{2r_{yz}}\sqn{x} + I_{\cL}(x)
	\\&\aeq{uses the definition of $F(x)$ in \cref{eq:FG}}
	\min_{x \in (\sX)^n} \sum_{i=1}^{n} \left(f_i(x_i) + \frac{r_x + 1/r_{yz}}{2}\sqn{x_i}\right) + I_{\cL}(x)
	\\&\aeq{uses \cref{eq:r_xyz}}
	\min_{x \in (\sX)^n} \sum_{i=1}^{n} \left(f_i(x_i) + \frac{r}{2}\sqn{x_i}\right) + I_{\cL}(x)
	\\&\aeq{uses the definition of $p(x)$ in \cref{eq:main} and the definition of $I_\cL(x)$}
	n \cdot \min_{x \in \sX} p(x).
\end{align*}
where \annotate.\qed

\newpage

\section{Proof of Theorems~\ref{thm:lb_sc} and~\ref{thm:lb}}\label{sec:lb_proof}

\subsection{The Hard Instance of Problem~(\ref{eq:main})}\label{sec:hard}

\paragraph{Compute nodes.}
In this proof, we consider the case when $\chi \geq 3$. The case $\chi < 3$ can be proven using the fixed-network argument of \citet{scaman2018optimal}.
We choose $n = 3\floor{\chi / 3}$, which implies that $n\geq 3$ and $n \bmod 3 = 0$. We also divide the set of nodes $\cV = \rng{1}{n}$ into the following three disjoint subsets: $\cV_1 = \rng{1}{n/3}$, $\cV_2= \rng{n/3+1}{2n/3}$ and $\cV_3 = \rng{2n/3+1}{n}$.

\paragraph{Objective functions.}
We fix an arbitrary odd integer $d \in \{3,5,7,\ldots\}$ and define functions $f_1(x), \ldots, f_n(x) \colon \sX \to \R$ as follows:
\begin{equation}\label{eq:lb_f}
	f_i(x) = \begin{cases}
		a{\textstyle \sum_{j=1}^{(d-1)/2}} h_{2j-1}(x) - a\<x,\basis{d}{1}> & i \in \cV_1 \\
		a{\textstyle \sum_{j=1}^{(d-1)/2}} h_{2j}(x)                        & i \in \cV_2 \\
		0                                                                   & i \in \cV_3 \\
	\end{cases},
\end{equation}
where $a > 0$ is an arbitrary constant and functions $h_1(x),\ldots,h_{d-1}(x)\colon \sX \to \R$ are defined as follows:
\begin{equation}\label{eq:lb_h}
	h_{j}(x) = \abs*{\<x,\basis{d}{j+1} - \basis{d}{j}>}.
\end{equation}
Consequently, the objective function $p(x)$ in problem~\eqref{eq:main} is given as follows:
\begin{equation}\label{eq:lb_p}
	p(x) = \frac{a}{3} \sum_{j=1}^{d-1} h_j(x) - \frac{a}{3}\<\basis{d}{1},x> + \frac{r}{2}\sqn{x}.
\end{equation}
We also define the subgradient oracles $\subgrad f_1(x),\ldots,\subgrad f_n(x) \colon \sX\to\sX$ as follows:
\begin{equation}\label{eq:subgrad_f}
	\subgrad f_i(x) = \begin{cases}
		a{\textstyle \sum_{j=1}^{(d-1)/2}} \subgrad h_{2j-1}(x) - a\basis{d}{1} & i \in \cV_1 \\
		a{\textstyle \sum_{j=1}^{(d-1)/2}} \subgrad h_{2j}(x)                   & i \in \cV_2 \\
		0                                                                       & i \in \cV_3 \\
	\end{cases},
\end{equation}
where $\subgrad h_1(x),\ldots,\subgrad h_{d-1}(x) \colon \sX\to\sX$ are the subgradient oracles associated with functions $h_1(x),\ldots,h_{d-1}(x)$, defined as follows:
\begin{equation}\label{eq:subgrad_h}
	\subgrad h_j(x) = \begin{cases}
		\basis{d}{j+1} - \basis{d}{j}   & \<\basis{d}{j+1},x> > \<\basis{d}{j},x> \\
		0                               & \<\basis{d}{j+1},x> = \<\basis{d}{j},x> \\
		\basis{d}{j}   - \basis{d}{j+1} & \<\basis{d}{j+1},x> < \<\basis{d}{j},x>
	\end{cases}.
\end{equation}

\paragraph{Time-varying network.}

We choose the time-varying network $\cG(\texe) = (\cV, \cE(\texe))$ to be a star-topology undirected graph with the time-varying center node $i_c(\texe) \in \cV$. Formally, we define the edges of the time-varying network $\cE(\texe) \subset \cV \times \cV$ as follows:
\begin{equation}\label{eq:E}
	\cE(\texe) = \bigcup_{i \in \cV, i \neq i_c(\texe)} \{(i,i_c(\texe)), (i_c(\texe), i)\}.
\end{equation}
We also specify the center node $i_c(\texe)$ at a given time $\texe \geq 0$ as follows:
\begin{equation}\label{eq:center}
	i_c(\texe) = 2n/3 + 1 + \left(\floor{{\texe}/{\tcom}} \bmod n/3\right).
\end{equation}
We choose the time-varying gossip matrix $\mW(\texe) \in \R^{n\times n}$ to be the Laplacian matrix of the graph $\cG(\texe)$. Formally, $\mW(\texe)$ is defined as follows:
\begin{equation}
	\mW(\texe)_{ij} = \frac{1}{n}\begin{cases}
		0             & i\neq j \text{ and } (i,j) \notin \cE(\texe) \\
		-1            & i\neq j \text{ and } (i,j) \in \cE(\texe)    \\
		\deg_i(\texe) & i = j
	\end{cases},
\end{equation}
where $\deg_i(\texe)$ denotes the degree of the node $i \in \cV$ in the graph $\cG(\texe)$, i.e.,
\begin{equation}
	\deg_i(\texe) = \abs*{\left\{j : (i,j) \in \cE(\texe) \right\}}.
\end{equation}
One can observe, that the time-varying gossip matrix $\mW(\texe)$ satisfies \Cref{ass:gossip}, in particular, $\ker \mW(\texe) = \ker \mW(\texe)^\top = \spanset(\{\ones_n\})$. Moreover, one can show that $\mW(\texe)$ is a symmetric matrix, and $\lmax(\mW(\texe)) = 1$ and $\lminp(\mW(\texe))=1/n \geq 1/\chi$. Hence, $\mW(\texe)$ satisfies \Cref{ass:chi}.

\subsection{Auxiliary Lemmas}\label{sec:lb_lem}

Further, we define linear spaces $\cK_{0},\ldots,\cK_d \subset \sX$ as follows:
\begin{equation}\label{eq:K}
	\cK_0 = \{0\}
	\quad\text{and}\quad
	\cK_j  = \spanset\left(\{\basis{d}{1},\ldots,\basis{d}{j}\}\right)
	\quad\text{for}\quad
	j \in \rng{1}{d}.
\end{equation}
In order to prove \Cref{thm:lb_sc,thm:lb}, we will use the following auxiliary lemmas. The proofs of these lemmas can be found in \Cref{sec:lb_lem_proof}. Furthermore, the proof of \Cref{thm:lb_sc} is contained in \Cref{sec:proof:lb_sc}, and the proof of \Cref{thm:lb} is contained in \Cref{sec:proof:lb}.

\begin{lemma}\label{lem:subgrad_span}
	For all $\texe \geq 0$, the following statements hold:
	\begin{enumerate}
		\item[\bf (i)] Let $i \in \cV_1$. Then, for all $j \in \rng{1}{(d-1)/2}$,
		      \begin{equation}
			      \mem{i}(\texe) \subset \cK_{2j}\quad\text{implies}\quad\memloc{i}(\texe+\tloc) \subset \cK_{2j}.
		      \end{equation}
		\item[\bf (ii)] Let $i \in \cV_2$. Then, for all $j \in \rng{0}{(d-1)/2}$,
		      \begin{equation}
			      \mem{i}(\texe) \subset \cK_{2j+1}\quad\text{implies}\quad\memloc{i}(\texe+\tloc) \subset \cK_{2j+1}.
		      \end{equation}
		\item[\bf (ii)] Let $i \in \cV_3$. Then, for all $j \in \rng{0}{d}$,
		      \begin{equation}
			      \mem{i}(\texe) \subset \cK_{j}\quad\text{implies}\quad\memloc{i}(\texe+\tloc) \subset \cK_{j}.
		      \end{equation}
	\end{enumerate}
\end{lemma}
The proof of \Cref{lem:subgrad_span} is contained in \Cref{sec:proof:subgrad_span}.

\begin{lemma}\label{lem:comm_span}
	Let $k \in \rng{0}{n(d-1)/6 - 1}$. Then, for all $\texe < (k+1)\tcom$, the following inclusion holds:
	\begin{equation}\label{eq:induction}
		\mem{i}(\texe) \subset \begin{cases}
			\cK_{2p + 2} & i \in \cV_1 \text{ or } \left(i \in \cV_3 \text{ and } i \leq 2n/3 + q + 1\right) \\
			\cK_{2p + 1} & i \in \cV_2 \text{ or } \left(i \in \cV_3 \text{ and } i > 2n/3 + q + 1\right)
		\end{cases},
	\end{equation}
	where $p = \floor{3k/n}$ and $q = k \bmod (n/3)$.
\end{lemma}
The proof of \Cref{lem:comm_span} is contained in \Cref{sec:proof:comm_span}.

\begin{lemma}\label{lem:lb_sol}
	Let functions $f_1,\ldots,f_n(x)$ be defined by \cref{eq:lb_f}. Then problem~\cref{eq:main} has a unique solution $x^* \in \sX$, which is given as follows:
	\begin{equation}\label{eq:lb_sol}
		x^* = \frac{a}{3rd} \ones_d.
	\end{equation}
	Moreover, for all $x \in \cK_{d-1}$, the following inequality holds:
	\begin{equation}\label{eq:lb_gap}
		p(x) - p(x^*) \geq \frac{a^2}{18rd}.
	\end{equation}
\end{lemma}
The proof of \Cref{lem:lb_sol} is contained in \Cref{sec:proof:lb_sol}.

\subsection{Proof of Theorem~\ref{thm:lb_sc}}\label{sec:proof:lb_sc}

\paragraph{Decentralized communication.}
\Cref{lem:comm_span} implies that $\mem{i}(\texe) \subset \cK_{d-1}$ as long as $\texe <  \tcom \cdot n(d-1)/6$. Hence, \Cref{lem:lb_sol} implies \cref{eq:lb_gap}
for all $x \in \mem{i}(\texe)$ as long as $\texe <  \tcom \cdot n(d-1)/6$. Let the constant $a > 0$ be chosen as follows:
\begin{equation}\label{eq:a}
	a = \frac{M}{2\sqrt{d}}.
\end{equation}
Then, each function $f_i(x)$ defined by \cref{eq:lb_f} is $M$-Lipschitz. Indeed, the case $i \in \cV_3$ is trivial. In the case when $i \in \cV_1$, we can prove the $M$-Lipschitz continuity of $f_i(x)$ as follows:
\begin{align*}
	f_i(x) - f_i(x')
	 & =
	a\smashoperator{\sum_{j=1}^{(d-1)/2}} \left(\abs*{\<x,\basis{d}{2j} - \basis{d}{2j-1}>} - \abs*{\<x',\basis{d}{2j} - \basis{d}{2j-1}>}\right) - a\<x-x',\basis{d}{1}>
	\\&\leq
	a\smashoperator{\sum_{j=1}^{(d-1)/2}} \abs*{\<x-x',\basis{d}{2j} - \basis{d}{2j-1}>} + a\abs*{\<x-x',\basis{d}{1}>}
	\\&\leq
	a\smashoperator{\sum_{j=1}^{(d-1)/2}} \left(
	\abs*{\<x-x',\basis{d}{2j} >}
	+
	\abs*{\<x-x',\basis{d}{2j-1}>}
	\right)
	+ a\abs*{\<x-x',\basis{d}{1}>}
	\\&=
	a\smashoperator{\sum_{j=1}^{d-1}}
	\abs*{\<x-x',\basis{d}{j} >}
	+ a\abs*{\<x-x',\basis{d}{1}>}
	\\&\leq
	2a\smashoperator{\sum_{j=1}^{d}}
	\abs*{\<x-x',\basis{d}{j} >}
	\leq
	2a\sqrt{d}\norm{x-x'}
	\leq
	M\norm{x-x'}.
\end{align*}
In the case when $i \in \cV_2$, we can prove the $M$-Lipschitz continuity of $f_i(x)$ similarly.

Without loss of generality, we assume $\epsilon \leq {M^2}/({576r})$ and define $d \in \{3,5,\ldots\}$ as follows:
\begin{equation}\label{eq:d}
	d = 2\floor*{\frac{M}{12\sqrt{r\epsilon}}} - 1.
\end{equation}
Using \cref{eq:a,eq:d}, for all $\texe <  \tcom \cdot n(d-1)/6$ and $x \in \mem{i}(\texe)$, we obtain
\begin{align*}
	p(x) - p(x^*) \geq \frac{M^2}{36rd^2} > \epsilon.
\end{align*}
Hence, to reach precision $p(x) - p(x^*) \leq \epsilon$ for some $x \in \mem{i}(\texe)$, it is necessary that $\texe$ satisfies
\begin{equation}\label{eq:lower_comm}
	\begin{split}
		\texe & \geq \tcom \cdot \frac{n(d-1)}{6}
		\\&=
		\tcom \cdot \floor*{\frac{\chi}{3}}\left(\floor*{\frac{M}{12\sqrt{r\epsilon}}} - 1\right)
		\\&\geq
		\tcom \cdot \frac{\chi}{3}\left(\frac{M}{12\sqrt{r\epsilon}} - 1\right)
		\\&=
		\Omega\left(\tcom \cdot \frac{M\chi}{\sqrt{r\epsilon}}\right).
	\end{split}
\end{equation}

\paragraph{Subgradient computation.}
We also need to prove that to reach precision $p(x) - p(x^*) \leq \epsilon$ for some $x \in \mem{i}(\texe)$, it is necessary that $\texe$ satisfies
\begin{equation}\label{eq:lower_subgrad}
	\texe \geq \Omega\left(\tloc \cdot \frac{M^2}{r\epsilon}\right).
\end{equation}
We can do this by providing an extended version of our hard problem instance, described in \Cref{sec:hard}. In particular, we consider the following instance of problem~\eqref{eq:main}:
\begin{equation}
	\min_{(x,x') \in \sX \times \R^{d'}} \frac{1}{n}\sum_{i=1}^{n} (f_i(x) + f_i'(x')) + \frac{r}{2}\sqn{x} + \frac{r}{2}\sqn{x'},
\end{equation}
where functions $f_1(x),\ldots,f_n(x) \colon \sX \to \R$ are defined in \Cref{sec:hard} by \cref{eq:lb_f}, and functions $f_1'(x'),\ldots,f_n'(x') \colon \R^{d'} \to \R$ are defined as follows:
\begin{equation}
	f_i'(x') = b\max_{j \in \rng{1}{d'}} \<\basis{d'}{j},x'>,
\end{equation}
where $b > 0$ is some constant. Then, by choosing an appropriate subgradient oracle $\subgrad f_i'(x')$ associated with each function $f'_i(x')$ (see Section~3.2.1~of \citet{nesterov2013introductory}) we can obtain both lower bounds~\eqref{eq:lower_comm} and~\eqref{eq:lower_subgrad}, which concludes the proof.\qed

\newpage

\subsection{Proof of Theorem~\ref{thm:lb}}\label{sec:proof:lb}

Our proof of \Cref{thm:lb} is very similar to the proof of \Cref{thm:lb_sc} with the following differences.
Let function $h_\delta(x)\colon\sX\to \R$ be the Huber function, which is defined as follows:
\begin{equation}
	h_\delta(x) = \sum_{j=1}^{d} h_\delta^j(\<\basis{d}{j},x>),\quad
	\text{where}\quad
	h_\delta^j(t) = \begin{cases}
		\frac{1}{2}t^2                       & \abs{t} \leq \delta \\
		\delta \abs{t} - \frac{1}{2}\delta^2 & \abs{t} > \delta
	\end{cases}.
\end{equation}
Note that function $h_\delta(x)$ is continuously differentiable and $(\sqrt{d}\delta)$-Lipschitz continuous.

In the proof of \Cref{thm:alg_sc} we used functions $f_1(x),\ldots,f_n(x)$ defined in \cref{eq:lb_f} of \Cref{sec:hard}. Here we use a slightly different choice, that is, functions $f_1(x),\ldots,f_n(x)$ are defined as follows:
\begin{equation}\label{eq:lb_f2}
	f_i(x) = h_\delta(x) + \begin{cases}
		a{\textstyle \sum_{j=1}^{(d-1)/2}} h_{2j-1}(x) - a\<x,\basis{d}{1}> & i \in \cV_1 \\
		a{\textstyle \sum_{j=1}^{(d-1)/2}} h_{2j}(x)                        & i \in \cV_2 \\
		0                                                                   & i \in \cV_3 \\
	\end{cases}.
\end{equation}
Consequently, our hard instance of problem~\eqref{eq:main}, which is described in \Cref{sec:hard}, turns into the following:
\begin{equation}\label{eq:main_lb}
	\min_{x \in \sX} \left[p(x) = \frac{a}{3} \sum_{j=1}^{d-1} h_j(x) - \frac{a}{3}\<\basis{d}{1},x> + c h_\delta(x)\right],
\end{equation}
where $c > 0$ is some constant, and functions $h_1(x),\ldots,h_{d-1}(x)$ are defined in \cref{eq:lb_h}.

One can show that \Cref{lem:subgrad_span,lem:comm_span} still hold true. We can also replace \Cref{lem:lb_sol} with the following \Cref{lem:lb_sol2}. The proof of this lemma is a trivial extension of the proof of \Cref{lem:lb_sol}, which uses the fact that $\nabla(\frac{1}{2}\sqn{\cdot})(x^*) = \nabla h_\delta(x^*)$ as long as $\delta$ and $x^*$ are defined by \cref{eq:lb_delta} and \cref{eq:lb_sol2}, respectively.
\begin{lemma}\label{lem:lb_sol2}
	Let $\delta$ be defined as follows:
	\begin{equation}\label{eq:lb_delta}
		\delta = \frac{a}{3cd}.
	\end{equation}
	Problem~\cref{eq:main_lb} has a solution $x^* \in \sX$, which is given as follows:
	\begin{equation}\label{eq:lb_sol2}
		x^* = \frac{a}{3cd} \ones_d.
	\end{equation}
	Moreover, for all $x \in \cK_{d-1}$, the following inequality holds:
	\begin{equation}\label{eq:lb_gap2}
		p(x) - p(x^*) \geq \frac{a^2}{18cd}.
	\end{equation}
\end{lemma}

One can also show that each function $f_i(x)$ defined in \cref{eq:lb_f2} is $M_f$-Lipschitz continuous, where $M_f$ is defined as follows:
\begin{equation}
	M_f = 2a\sqrt{d} + c\delta\sqrt{d} = 2a\sqrt{d} + {a}/({3\sqrt{d}}) \leq 3a\sqrt{d}.
\end{equation}
Let us choose $a$ and $c$ as follows:
\begin{equation}
	a = \frac{M}{3\sqrt{d}} \quad\text{and}\quad c = \frac{M}{9Rd}.
\end{equation}
This choice of $a$ and c implies $M_f \leq M$ and $\norm{x^*} \leq R$. Moreover, \cref{eq:lb_gap2} implies
\begin{equation}
	p(x) - p(x^*) \geq  \frac{MR}{18d}
\end{equation}
as long as $x^* \in \cK_{d-1}$. Next, without loss of generality we can assume $\epsilon \leq (MR)/72$ and choose $d \in \{3,5,\ldots\}$ as follows:
\begin{equation}
	d  = 2\floor*{\frac{MR}{36\epsilon}} - 1,
\end{equation}
which, for all $x \in \mem{i}(\texe)$, implies
\begin{align*}
	p(x) - p(x^*) > \epsilon
\end{align*}
as long as $\texe$ satisfies
\begin{equation}
	\begin{split}
		\texe \geq \tcom \cdot \frac{n(d-1)}{6}=
		\Omega\left(\tcom \cdot \frac{MR\chi}{\epsilon}\right),
	\end{split}
\end{equation}
which concludes the proof.\qed

\newpage

\section{Proofs of Lemmas from Section~\ref{sec:lb_lem}}\label{sec:lb_lem_proof}

\subsection{Proof of Lemma~\ref{lem:subgrad_span}}\label{sec:proof:subgrad_span}

\paragraph{Statement (i).}
Let $i \in \cV_1$ and $x \subset \cK_{2j}$ for $j \in \rng{1}{(d-1)/2}$.
Then for $l \geq 2j+1$ we obtain $\<\basis{d}{l+1} - \basis{d}{l},x> = 0$, which implies $\subgrad h_l(x) = 0$ due to \cref{eq:subgrad_h}. Hence, we obtain the following:
\begin{align*}
	\tfrac{1}{a}\subgrad f_i(x)
	 & \aeq{uses \cref{eq:subgrad_f}}
	\subgrad h_1(x) + \subgrad h_3(x) + \cdots + \subgrad h_{d-2} (x) - \basis{d}{1}
	\\&\aeq{uses the fact that $\subgrad h_l(x) = 0$ for $l \geq 2j+1$}
	\subgrad h_1(x) + \subgrad h_3(x) + \cdots + \subgrad h_{2j-1} (x) - \basis{d}{1}
	\\&\asubset{uses \cref{eq:subgrad_h}} \spanset\left(\{\basis{d}{1},\basis{d}{2}\}\cup \cdots\cup \{\basis{d}{2j-1},\basis{d}{2j}\} \right)
	\\&\asubset{uses the definition of $\cK_{2j}$ in \cref{eq:K}}
	\cK_{2j},
\end{align*}
where \annotate. Hence, $\mem{i}(\texe) \subset \cK_{2j}$ implies $\memloc{i}(\texe+\tloc) \subset \cK_{2j}$ by the definition of $\memloc{i}(\cdot)$ in \cref{eq:memloc}.

\paragraph{Statement (ii).}
Let $i \in \cV_2$ and $x \subset \cK_{2j+1}$ for $j \in \rng{0}{(d-1)/2}$.
Then for $l \geq 2j+2$ we obtain $\<\basis{d}{l+1} - \basis{d}{l},x> = 0$, which implies $\subgrad h_l(x) = 0$ due to \cref{eq:subgrad_h}. Hence, we obtain the following:
\begin{align*}
	\tfrac{1}{a}\subgrad f_i(x)
	 & \aeq{uses \cref{eq:subgrad_f}}
	\subgrad h_2(x) + \subgrad h_4(x) + \cdots + \subgrad h_{d-1} (x)
	\\&\aeq{uses the fact that $\subgrad h_l(x) = 0$ for $l \geq 2j+2$}
	\subgrad h_2(x) + \subgrad h_4(x) + \cdots + \subgrad h_{2j} (x)
	\\&\asubset{uses \cref{eq:subgrad_h}} \spanset\left(\{\basis{d}{2},\basis{d}{3}\}\cup \cdots\cup \{\basis{d}{2j},\basis{d}{2j+1}\} \right)
	\\&\asubset{uses the definition of $\cK_{2j+1}$ in \cref{eq:K}}
	\cK_{2j},
\end{align*}
where \annotate. Hence, $\mem{i}(\texe) \subset \cK_{2j+1}$ implies $\memloc{i}(\texe+\tloc)\subset \cK_{2j+1}$ by the definition of $\memloc{i}(\cdot)$ in \cref{eq:memloc}.

\paragraph{Statement (iii).} This statement is trivially implied by the definition of $\subgrad f_i(x)$ in \cref{eq:subgrad_f} and the definition of $\memloc{i}(\cdot)$ in \cref{eq:memloc}. \qed

\subsection{Proof of Lemma~\ref{lem:comm_span}}\label{sec:proof:comm_span}

We prove the lemma using the induction on $k$.

\paragraph{Base case: $k=0$.} In this case, we assume $\texe < (k+1)\tcom = \tcom$. Hence, for all $i \in \cV$, we obtain $\memcom{i}(\texe) = \varnothing$ and $\mem{i}(\texe) \subset \memloc{i}(\texe)$. Using \Cref{lem:subgrad_span} and the fact that $\memcom{i}(\texe) = \varnothing$, we can easily obtain
\begin{align*}
	\mem{i}(\texe) \subset \memloc{i}(\texe) \subset \begin{cases}
		                                                 \cK_2 & i \in \cV_1 \\
		                                                 \cK_1 & i \in \cV_2 \\
		                                                 \cK_0 & i \in \cV_3 \\
	                                                 \end{cases},
\end{align*}
which implies the desired \cref{eq:induction} for $k=p=q=0$.

\paragraph{Induction hypothesis.} Let $k' \in \{0,1,2,\ldots\}$. We assume that \cref{eq:induction} holds for all $\texe < (k'+1)\tcom$, that is,
\begin{equation}\label{eq:induction_hypothesis}
	\mem{i}(\texe) \subset
	\begin{cases}
		\cK_{2p' + 2} & i \in \cV_1 \text{ or } \left(i \in \cV_3 \text{ and } i \leq 2n/3 + q' + 1\right) \\
		\cK_{2p' + 1} & i \in \cV_2 \text{ or } \left(i \in \cV_3 \text{ and } i > 2n/3 + q' + 1\right)
	\end{cases},
\end{equation}
where $p' = \floor{3k'/n}$ and $q' = k' \bmod (n/3)$.

\paragraph{Induction step.} We assume that the induction hypothesis~\eqref{eq:induction_hypothesis} is true. Our goal is to prove that \cref{eq:induction} holds for $k = k'+ 1$.
When $0 \leq \texe < k\tcom$, the desired \cref{eq:induction} is implied by the induction hypothesis~\eqref{eq:induction_hypothesis}. Thus, we can assume $k\tcom\leq \texe < (k+1)\tcom$.
Further, we consider two cases: $q \neq 0$ and $q = 0$.

\paragraph{Induction step, case $q \neq 0$.} In this case, $p = p'$ and $q = q'+1$.

	{\bf Part  (i).} First, we consider the case
\begin{equation}\label{eq:aux6}
	k\tcom\leq \texe <  \min\{(k+1)\tcom,k\tcom + \tloc\}.
\end{equation}
\Cref{eq:aux6} implies $\texe - \tloc < (k'+1)\tcom$ and $\texe - \tcom < (k'+1)\tcom$. Using the induction hypothesis~\eqref{eq:induction_hypothesis} and the fact that $p' = p$ and $q' = q-1$, we get
\begin{equation}\label{eq:aux7}
	\mem{i}(\texe - \tloc),\;\mem{i}(\texe - \tcom)
	\subset
	\begin{cases}
		\cK_{2p + 2} & i \in \cV_1 \text{ or } \left(i \in \cV_3 \text{ and } i \leq 2n/3 + q\right) \\
		\cK_{2p + 1} & i \in \cV_2 \text{ or } \left(i \in \cV_3 \text{ and } i > 2n/3 + q\right)
	\end{cases}.
\end{equation}
Hence, using \Cref{lem:subgrad_span}, we obtain
\begin{equation}\label{eq:aux8}
	\memloc{i}(\texe)
	\subset
	\begin{cases}
		\cK_{2p + 2} & i \in \cV_1 \text{ or } \left(i \in \cV_3 \text{ and } i \leq 2n/3 + q\right) \\
		\cK_{2p + 1} & i \in \cV_2 \text{ or } \left(i \in \cV_3 \text{ and } i > 2n/3 + q\right)
	\end{cases}.
\end{equation}

\Cref{eq:aux6,eq:center} imply $i_c(\texe) = 2n/3 + q + 1$. Hence, using \cref{eq:aux7}, we get
\begin{equation*}
	\mem{i_c(\texe)}(\texe - \tcom)\subset\cK_{2p + 1}.
\end{equation*}
For $i \neq i_c(\texe)$, using \cref{eq:aux7,eq:memcom}, we get
\begin{equation}\label{eq:aux9}
	\memcom{i}(\texe) = \spanset\left(\mem{i_c(\texe)}(\texe-\tcom)\right) \subset \cK_{2p + 1}.
\end{equation}
For $i = i_c(\texe) = 2n/3 + q + 1$, using \cref{eq:aux7,eq:memcom}, we get
\begin{equation}\label{eq:aux10}
	\memcom{i_c(\texe)}(\texe) = \spanset\left(\bigcup_{j \neq i_c(\texe)}
	\mem{j}(\texe-\tcom)\right) \subset \cK_{2p + 2}.
\end{equation}
Hence, using \cref{eq:aux9,eq:aux10}, for all $i \in \cV$, we obtain
\begin{equation}\label{eq:aux11}
	\memcom{i}(\texe)
	\subset
	\begin{cases}
		\cK_{2p + 2} & i  = 2n/3 + q + 1     \\
		\cK_{2p + 1} & i  \neq  2n/3 + q + 1
	\end{cases}.
\end{equation}
Now, we combine \cref{eq:aux8,eq:aux11}, and obtain
\begin{align*}
	\mem{i}(\texe) & \subset \memloc{i}(\texe) \cup \memcom{i}(\texe)
	\subset
	\begin{cases}
		\cK_{2p + 2} & i \in \cV_1 \text{ or } \left(i \in \cV_3 \text{ and } i \leq 2n/3 + q + 1\right) \\
		\cK_{2p + 1} & i \in \cV_2 \text{ or } \left(i \in \cV_3 \text{ and } i > 2n/3 + q + 1\right)
	\end{cases}.
\end{align*}
Thus, we were able to prove \cref{eq:induction} for $\texe$ satisfying \eqref{eq:aux6}.

{\bf Part (ii).} We can prove the general case
\begin{equation*}
	k\tcom\leq \texe <  \min\{(k+1)\tcom,k\tcom + l\tloc\}
\end{equation*}
for arbitrary $l \in \{1,2,\ldots\}$ using the induction on $l$. The only difference compared to the proof in the previous part is in \cref{eq:aux7}, which will change to
\begin{align*}
	\mem{i}(\texe - \tloc)
	\subset
	\begin{cases}
		\cK_{2p + 2} & i \in \cV_1 \text{ or } \left(i \in \cV_3 \text{ and } i \leq 2n/3 + q + 1\right) \\
		\cK_{2p + 1} & i \in \cV_2 \text{ or } \left(i \in \cV_3 \text{ and } i > 2n/3 + q + 1\right)
	\end{cases}.
\end{align*}
However, \cref{eq:aux8} will not change due to \Cref{lem:subgrad_span}. Hence, the rest of the proof will also remain unchanged.

\paragraph{Induction step, case $q = 0$.} In this case $p = p'+1$ and $q' = n/3 - 1$.

	{\bf Part (i).}
First, we consider the case
\begin{equation}\label{eq:aux12}
	k\tcom\leq \texe <  \min\{(k+1)\tcom,k\tcom + \tloc\}.
\end{equation}
\Cref{eq:aux12} implies $\texe - \tloc < (k'+1)\tcom$ and $\texe - \tcom < (k'+1)\tcom$. Using the induction hypothesis~\eqref{eq:induction_hypothesis} and the fact that $p' = p-1$ and $q' = n/3-1$, we get
\begin{equation}\label{eq:aux13}
	\mem{i}(\texe - \tloc),\;\mem{i}(\texe - \tcom)
	\subset
	\begin{cases}
		\cK_{2p}     & i \in \cV_1 \text{ or } i \in \cV_3 \\
		\cK_{2p - 1} & i \in \cV_2
	\end{cases}.
\end{equation}

\Cref{eq:aux12,eq:center} imply $i_c(\texe) = 2n/3 + 1$. Using \cref{eq:aux13}, we get
\begin{equation*}
	\mem{i_c(\texe)}(\texe - \tcom)\subset\cK_{2p}.
\end{equation*}
For $i \neq i_c(\texe)$, using \cref{eq:aux13,eq:memcom}, we get
\begin{equation}\label{eq:aux14}
	\memcom{i}(\texe) = \spanset\left(\mem{i_c(\texe)}(\texe-\tcom)\right) \subset \cK_{2p}.
\end{equation}
For $i = i_c(\texe) = 2n/3 + 1$, using \cref{eq:aux13,eq:memcom}, we get
\begin{equation}\label{eq:aux15}
	\memcom{i_c(\texe)}(\texe) = \spanset\left(\bigcup_{j \neq i_c(\texe)}
	\mem{j}(\texe-\tcom)\right) \subset \cK_{2p}.
\end{equation}
Hence, using \cref{eq:aux14,eq:aux15}, for all $i \in \cV$, we obtain
\begin{equation}\label{eq:aux16}
	\memcom{i}(\texe)
	\subset
	\cK_{2p}.
\end{equation}
Using \Cref{lem:subgrad_span}, from \cref{eq:aux13} we obtain
\begin{equation}\label{eq:aux17}
	\memloc{i}(\texe)
	\subset
	\begin{cases}
		\cK_{2p}     & i \in \cV_1 \text{ or } i \in \cV_3 \\
		\cK_{2p - 1} & i \in \cV_2
	\end{cases}.
\end{equation}
Hence, using \cref{eq:aux16,eq:aux17}, for all $i \in \cV$, we obtain
\begin{equation}\label{eq:aux18}
	\mem{i}(\texe) \subset \memloc{i}(\texe) \cup \memcom{i}(\texe)
	\subset
	\cK_{2p},
\end{equation}
which implies \cref{eq:induction} for $\texe$ satisfying \eqref{eq:aux12}.

{\bf Part (ii).}
Next, we consider the case
\begin{equation}\label{eq:aux19}
	k\tcom + \tloc \leq \texe < \min\{(k+1)\tcom,k\tcom + 2\tloc\}.
\end{equation}
\Cref{eq:aux16} still holds for all $i \in \cV$ and $\texe$ satisfying \cref{eq:aux19}. From \cref{eq:aux19,eq:aux18}, for all $i \in \cV$, we obtain
\begin{equation*}
	\mem{i}(\texe-\tloc) \subset \cK_{2p},
\end{equation*}
which, due to \Cref{lem:subgrad_span}, implies the following:
\begin{equation}\label{eq:aux20}
	\memloc{i}(\texe) \subset \begin{cases}
		\cK_{2p}   & i \in \cV_1 \text{ or } i \in \cV_3 \\
		\cK_{2p+1} & i \in \cV_2
	\end{cases}.
\end{equation}
Hence, using \cref{eq:aux16,eq:aux20}, we obtain
\begin{equation}\label{eq:aux21}
	\mem{i}(\texe) \subset \memloc{i}(\texe) \cup \memcom{i}(\texe)
	\subset
	\begin{cases}
		\cK_{2p}   & i \in \cV_1 \text{ or } i \in \cV_3 \\
		\cK_{2p+1} & i \in \cV_2
	\end{cases},
\end{equation}
which implies \cref{eq:induction} for $\texe$ satisfying \eqref{eq:aux19}.

{\bf Part(iii).} We can prove the general case
\begin{equation}\label{eq:aux22}
	k\tcom + l\tloc \leq \texe < \min\{(k+1)\tcom,k\tcom + (l+1)\tloc\}
\end{equation}
for $l \in \{2,3,\ldots\}$ using the induction on $l$. There will be no differences compared to the proof in the previous part. Indeed, \cref{eq:aux16,eq:aux21} will still hold for all $i \in \cV$ and $\texe$ satisfying \cref{eq:aux22}. \qed

\newpage

\subsection{Proof Lemma~\ref{lem:lb_sol}}\label{sec:proof:lb_sol}

One can show, that $x^*$ defined in \cref{eq:lb_sol} is indeed the unique minimizer of the function $p(x)$ defined in \cref{eq:lb_p}. Moreover, we can obtain the following:
\begin{align*}
	p(x^*) = -\frac{a^2}{18rd}.
\end{align*}
We can lower-bound function $p(x)$ as follows:
\begin{align*}
	p(x)
	 & =
	\frac{r}{2}\sqn{x} - \frac{a}{3}\<\basis{d}{1},x> + \frac{a}{3} \sum_{j=1}^{d-1} \abs*{\<x,\basis{d}{j+1} - \basis{d}{j}>}
	\\&\geq
	-\frac{a}{3}\abs*{\<\basis{d}{1},x>}
	+\frac{a}{3} \sum_{j=1}^{d-1} \left(\abs*{\<x,\basis{d}{j}>} - \abs*{\<x,\basis{d}{j+1}>}\right)
	\\&=
	-\frac{a}{3}\abs*{\<\basis{d}{d},x>}
	\\&=0
\end{align*}
as long as $x \in \cK_{d-1}$. Hence, for all $x \in \cK_j$, we obtain
\begin{align*}
	p(x) - p(x^*) \geq \frac{a^2}{18rd},
\end{align*}
which concludes the proof.\qed

\newpage
\section{Proof of Theorems~\ref{thm:alg_sc} and~\ref{thm:alg}}\label{sec:alg_proof}

\subsection{Auxiliary Lemmas}\label{sec:alg_lem}

In order to prove \Cref{thm:alg_sc,thm:alg}, we will use the following auxiliary lemmas. The proofs of these lemmas can be found in \Cref{sec:alg_lem_proof}. Furthermore, the proof of \Cref{thm:alg_sc} is contained in \Cref{sec:proof:alg_sc}, and the proof of \Cref{thm:alg} is contained in \Cref{sec:proof:alg}.

\begin{lemma}\label{lem:sol}
	Under \Cref{ass:cvx,,ass:lip,,ass:sol}, let $r > 0$ (strongly convex case). Then there exists a solution $(w^*,y^*,z^*) \in \cL \times (\sX)^n \times \cL^\perp$ to problem~\eqref{eq:main}, which satisfies the following conditions
	\begin{equation}
		0 \in \partial_x Q(w^*,y^*,z^*),\quad
		0 = \nabla_y Q(w^*,y^*,z^*),\quad
		\cL \ni \nabla_z Q(w^*,y^*,z^*).
	\end{equation}
	Moreover, the following inequalities hold:
	\begin{equation}
		\sqn{w^*} \leq {nM^2}/{r^2},\quad
		\sqn{y^*} \leq (1+r_x/r)^2 nM^2,\quad
		\sqn{z^*} \leq 4nM^2.
	\end{equation}
\end{lemma}
The proof of \Cref{lem:sol} is contained in \Cref{sec:proof:sol}.

\begin{lemma}\label{lem:inner}
	Under \Cref{ass:cvx,ass:lip}, let $\eta_x^0,\ldots,\eta_x^{K-1}$ and $\beta_0,\ldots,\beta_{K-1}$ be chosen as follows:
	\begin{equation}\label{eq:eta_x_beta_sigma_k}
		\eta_x^k = 1/(\tau_x^kT),\quad
		\beta_k = r_x,\quad
		\sigma_k = \tau_x^k / (2\tau_x^k + \beta_k) \quad
		\text{for}\quad
		k \in \rng{0}{K-1}.
	\end{equation}
	Then, for all $x \in (\sX)^n$ and $k \in \rng{0}{K-1}$, the following inequality holds:
	\begin{equation}
		\begin{aligned}
			(\tau_x^k + \tfrac{1}{2}r_x)\sqn{x^{k+1} - x}
			 & \leq
			\tau_x^k\sqn{x^k - x}
			+{2n M^2}/({\tau_x^k T})
			\\&-
			\left(
			F(\tilde{x}^{k+1}) - F(x)
			-\<y^{k+1},\tilde{x}^{k+1} - x>
			+\tfrac{1}{2}\tau_x^k\sqn{\tilde{x}^{k+1} - x^k}
			\right).
		\end{aligned}
	\end{equation}
\end{lemma}
The proof of \Cref{lem:inner} is contained in \Cref{sec:proof:inner}.

\begin{lemma}\label{lem:P}
	Under \Cref{ass:gossip},
	for all $k \in \rng{0}{K-1}$, the iterates of \Cref{alg} satisfy
	\begin{equation}
		\mP z^k = z^k,\quad
		\mP\ol{z}^{k+1} = \ol{z}^{k+1},\quad
		\mP\uln{z}^k = \uln{z}^k,
	\end{equation}
	where $\mP \in \R^{nd \times nd}$ is the orthogonal projection matrix onto $\cL^\perp$, which is given as follows:
	\begin{equation}\label{eq:P}
		\mP = (\mI_n - \tfrac{1}{n}\ones_n\ones_n^\top)\otimes \mI_d.
	\end{equation}
\end{lemma}
The proof of \Cref{lem:P} is contained in \Cref{sec:proof:P}.

\begin{lemma}\label{lem:EF}
	Under \Cref{ass:gossip,ass:chi}, for all $k \in \rng{0}{K-1}$ the following inequality holds:
	\begin{equation}
		\sqn{\eta_z^k m^k}_{\mP} \leq 2\chi\sqn{\eta_z^k m^k}_{\mP}
		-2\chi\sqn{\eta_z^{k+1} m^{k+1}}_{\mP}
		+4\chi^2\sqn{\eta_z^k g_z^k}_{\mP}.
	\end{equation}
\end{lemma}
The proof of \Cref{lem:EF} is contained in \Cref{sec:proof:EF}.

\begin{lemma}\label{lem:descent}
	Under \Cref{ass:gossip,ass:chi}, let parameters $\theta_z^0,\ldots,\theta_z^{K-1}$ be chosen as follows:
	\begin{equation}\label{eq:theta_z_k}
		\theta_z^k = 1/(2r_{yz})
		\quad\text{for}\quad
		k = 0,\ldots,K-1.
	\end{equation}
	Then, for all $k \in \rng{0}{K-1}$, the following inequality holds:
	\begin{equation}
		0 \leq -\alpha_k^{-1}\left(
		\<\ol{z}^{k+1} - \uln{z}^k,g_z^k>
		+r_{yz}\sqn{\ol{z}^{k+1} - \uln{z}^k}
		\right)
		-({4\alpha_k\chi r_{yz}})^{-1}\sqn{g_z^k}_{\mP}.
	\end{equation}
\end{lemma}
The proof of \Cref{lem:descent} is contained in \Cref{sec:proof:descent}.

\begin{lemma}\label{lem:outer}
	Under \Cref{ass:cvx,ass:lip} and under conditions of \Cref{lem:inner,lem:descent},
	let parameters $\alpha_0,\ldots,\alpha_{K-1}$ and $\gamma_0,\ldots,\gamma_{K-1}$ be chosen as follows:
	\begin{equation}\label{eq:alpha_gamma_k}
		\alpha_k = 3/(k+3),\quad
		\gamma_k = (k+2)/(k+3)\quad
		\text{for}\quad
		k = 0,\ldots,K-1.
	\end{equation}
	Let parameters $\tau_x^0,\ldots,\tau_x^{K-1}$, $\eta_y^0,\ldots,\eta_y^{K-1}$, and $\eta_z^0,\ldots,\eta_z^{K-1}$ be chosen as follows:
	\begin{equation}\label{eq:tau_x_eta_yz_k}
		\tau_x^k = \tau_x\alpha_k^{-1},\quad
		\eta_y^k = \eta_y\alpha_k^{-1},\quad
		\eta_z^k = \eta_z\alpha_k^{-1}\quad
		\text{for}\quad
		k = 0,\ldots,K-1,
	\end{equation}
	where $\tau_x$, $\eta_y$ and $\eta_z$ are defined as follows:
	\begin{equation}\label{eq:tau_x_eta_yz}
		\tau_x = \tfrac{1}{2}r_x,
		\quad
		\eta_y = (4r_{yz})^{-1},
		\quad
		\eta_z = (10r_{yz}\chi^2)^{-1},
		\quad
		r_x = \tfrac{2}{3}r,
		\quad
		r_{yz} = 3/r.
	\end{equation}
	Let parameters $\lambda_1,\ldots,\lambda_K$ be chosen as follows:
	\begin{equation}\label{eq:lambda_k}
		\lambda_K = \alpha_{K-1}^{-2}
		\quad\text{and}\quad
		\lambda_k = \alpha_{k-1}^{-2} + \alpha_k^{-1} - \alpha_k^{-2}
		\quad\text{for}\quad k = 1,\ldots,K-1.
	\end{equation}
	Let the input of \Cref{alg} be chosen as follows:
	\begin{equation}\label{eq:input}
		x^0 = 0,
		\quad
		y^0 = 0,
		\quad
		z^0 = 0,
		\quad
		m^0 = 0.
	\end{equation}
	Then, for all $x,y \in (\sX)^n$ and $z \in \cL^\perp$, the following inequality holds:
	\begin{equation}
		Q(x_a^{K},y,z) - Q(x,y_a^K,z_a^K)
		\leq
		\frac{2}{K^2}
		\left(r\sqn{x}
		+\frac{18}{r}\sqn{y}
		+\frac{45\chi^2}{r}\sqn{z}
		\right)
		+\frac{72n M^2}{r KT}.
	\end{equation}
\end{lemma}
The proof of \Cref{lem:outer} is contained in \Cref{sec:proof:outer}.

\subsection{Proof of Theorem~\ref{thm:alg_sc}}\label{sec:proof:alg_sc}

We can upper-bound $\frac{r_x}{2}\sqn{x_a^K - w^*}$, where $w^*$ is defined in \Cref{lem:sol}, as follows:
\begin{align*}
	\frac{r_x}{2}\sqn{x_a^K - w^*}
	 & \aleq{uses \Cref{lem:sol} and the strong convexity of $Q(x,y,z)$ in $x$}
	Q(x_a^{K},y^*,z^*) - Q(w^*,y^*,z^*)
	\\&\aleq{use \Cref{lem:sol}}
	Q(x_a^{K},y^*,z^*) - Q(w^*,y_a^K,z_a^K)
	\\&\aleq{uses \Cref{lem:outer}}
	\frac{2}{K^2}
	\left(r\sqn{w^*}
	+\frac{18}{r}\sqn{y^*}
	+\frac{45\chi^2}{r}\sqn{z^*}
	\right)
	+\frac{72n M^2}{r KT}
	\\&\aleq{use \Cref{lem:sol}}
	\frac{2}{K^2}
	\left(\frac{nM^2}{r}
	+\frac{18(1+r_x/r)^2 nM^2}{r}
	+\frac{180n\chi^2M^2}{r}
	\right)
	+\frac{72n M^2}{r KT}
\end{align*}
where \annotate. Using the definition of $r_x$ in \cref{eq:tau_x_eta_yz}
\begin{align*}
	r\sqn{x_a^K - w^*}
	 & \leq
	\frac{6}{K^2}
	\left(
	\frac{51 nM^2}{r}
	+\frac{180n\chi^2M^2}{r}
	\right)
	+\frac{72n M^2}{r KT}
	\\&\leq
	\frac{1386 n \chi^2 M^2}{rK^2}
	+\frac{72 n M^2}{r KT}.
\end{align*}

Next, we can upper-bound $n(p(x_o^K) - p(x^*))$ as follows:
\begin{align*}
	n(p(x_o^K) - p(x^*))
	 & \aeq{uses the definition of $p(x)$ in \cref{eq:main}}
	\sum_{i=1}^{n} \left(f_i(x_{o}^K) - f_i(x^*) + \frac{r}{2}\sqn{x_o^K} - \frac{r}{2}\sqn{x^*}\right)
	\\&\aeq{uses the definition of $x_o^K$ on \cref{line:output} of \Cref{alg}}
	\sum_{i=1}^{n} \left(f_i(x_{o}^K) - f_i(x^*) + \frac{r}{2}\sqn{{\textstyle \frac{1}{n}\sum_{j=1}^{n}x_{a,j}^K}} - \frac{r}{2}\sqn{x^*}\right)
	\\&\aleq{uses the convexity of $\sqn{\cdot}$}
	\sum_{i=1}^{n} \left(f_i(x_{o}^K) - f_i(x^*) + \frac{r}{2}\sqn{x_{a,i}^K} - \frac{r}{2}\sqn{x^*}\right)
	\\&\aeq{uses the definition of $w^*$ in \cref{eq:w_star}}
	\sum_{i=1}^{n} \left(f_i(x_{o}^K) - f_i(x^*) \right)+ \frac{r}{2}\sqn{x_{a}^K} - \frac{r}{2}\sqn{w^*}
	\\&\aleq{uses \Cref{ass:lip}}
	\sum_{i=1}^{n} \left(f_i(x_{a,i}^K) - f_i(x^*) + M\norm{x_{a,i}^K - x_{o}^K}\right)+ \frac{r}{2}\sqn{x_{a}^K} - \frac{r}{2}\sqn{w^*}
	\\&\aeq{uses the definition of function $F(x)$ in \cref{eq:FG} and \cref{eq:r_xyz}}
	F(x_a^K) - F(w^*) + \frac{1}{2r_{yz}}\sqn{x_{a}^K} - \frac{1}{2r_{yz}}\sqn{w^*}
	+\sum_{i=1}^{n} M\norm{x_{a,i}^K - x_{o}^K}
	\\&\aleq{uese the Cauchy-Schwarz inequality}
	F(x_a^K) - F(w^*) + \frac{1}{2r_{yz}}\sqn{x_{a}^K} - \frac{1}{2r_{yz}}\sqn{w^*}
	\\&
	+\sqrt{{\textstyle\sum_{i=1}^{n}} M^2}\sqrt{{\textstyle\sum_{i=1}^{n}} \sqn{x_{a,i}^K - x_{o}^K}}
	\\&\aeq{uses the definition of $\mP$ in \cref{eq:P}}
	F(x_a^K) - F(w^*) + \frac{1}{2r_{yz}}\sqn{x_{a}^K} - \frac{1}{2r_{yz}}\sqn{w^*}
	+\sqrt{n}M \norm{x_a^K}_{\mP}
\end{align*}
where \annotate.

Next, for arbitrary $z \in \cL^\perp$ we define $y = -r_{yz}^{-1} x_a^K - z$. Then, we get $\nabla_y Q(x_a^K,y,z) = 0$ and $Q(x_a^K,y,z) = F(x_a^K) + \frac{1}{2r_{yz}}\sqn{x_a^K} - \<x_a^K,z>$. Plugging this into the previous upper-bound gives the following:
\begin{align*}
	n(p(x_o^K) - p(x^*))
	 & \leq
	Q(x_a^K,y,z) - F(w^*) - \frac{1}{2r_{yz}}\sqn{w^*} + \<x_a^K,z> + \sqrt{n}M \norm{x_a^K}_{\mP}
	\\&\aeq{uses the definition of $y^*$ in \cref{eq:y_star} and the definition of $z^*$ in \cref{eq:z_star}}
	Q(x_a^K,y,z) - Q(w^*,y^*,z^*) + \<x_a^K,z> + \sqrt{n}M \norm{x_a^K}_{\mP}
	\\&\aleq{uses \Cref{lem:sol}}
	Q(x_a^K,y,z) - Q(w^*,y_a^K,z_a^K) + \<x_a^K,z> + \sqrt{n}M \norm{x_a^K}_{\mP}
\end{align*}
where \annotate.

Next, we choose $z \in \cL^\perp$ as follows:
\begin{equation}
	z = \begin{cases}
		-\sqrt{n}M\norm{\mP x_a^K}^{-1}\mP x_a^K & x_a^K \neq 0 \\
		0                                        & x_a^K = 0
	\end{cases}.
\end{equation}
Then, $\<x_a^K,z> = -\sqrt{n}M \norm{x_a^K}_{\mP}$ and we obtain the following:
\begin{align*}
	\mind{2em}
	n(p(x_o^K) - p(x^*))
	\\&\leq
	Q(x_a^K,y,z) - Q(w^*,y_a^K,z_a^K)
	\\&\aleq{uses \Cref{lem:outer}}
	\frac{2}{K^2}
	\left(r\sqn{w^*}
	+\frac{18}{r}\sqn{y}
	+\frac{45\chi^2}{r}\sqn{z}
	\right)
	+\frac{72n M^2}{r KT}
	\\&\aeq{uses our choice of $y$}
	\frac{2}{K^2}
	\left(r\sqn{w^*}
	+\frac{18}{r}\sqn{r_{yz}^{-1} x_a^K + z}
	+\frac{45\chi^2}{r}\sqn{z}
	\right)
	+\frac{72n M^2}{r KT}
	\\&=
	\frac{2}{K^2}
	\left(r\sqn{w^*}
	+\frac{18}{r}\sqn{r_{yz}^{-1} (x_a^K - w^* + w^*) + z}
	+\frac{45\chi^2}{r}\sqn{z}
	\right)
	+\frac{72n M^2}{r KT}
	\\&\aleq{uses the parallelogram rule and Young's inequality}
	\frac{2}{K^2}
	\left(r\sqn{w^*}
	+\frac{54}{rr_{yz}^2}\sqn{x_a^K - w^*}
	+\frac{54}{rr_{yz}^2}\sqn{w^*}
	+\frac{54}{r}\sqn{z}
	+\frac{45\chi^2}{r}\sqn{z}
	\right)
	+\frac{72n M^2}{r KT}
	\\&\leq
	\frac{2}{K^2}
	\left(r\sqn{w^*}
	+\frac{54}{rr_{yz}^2}\sqn{x_a^K - w^*}
	+\frac{54}{rr_{yz}^2}\sqn{w^*}
	+\frac{99\chi^2}{r}\sqn{z}
	\right)
	+\frac{72n M^2}{r KT}
	\\&\aleq{uses our choice of $z$}
	\frac{2}{K^2}
	\left(r\sqn{w^*}
	+\frac{54}{rr_{yz}^2}\sqn{x_a^K - w^*}
	+\frac{54}{rr_{yz}^2}\sqn{w^*}
	+\frac{99n\chi^2M^2}{r}
	\right)
	+\frac{72n M^2}{r KT}
	\\&\aleq{uses \Cref{lem:sol}}
	\frac{2}{K^2}
	\left(
	\frac{nM^2}{r}
	+\frac{54nM^2}{r^3r_{yz}^2}
	+\frac{54}{rr_{yz}^2}\sqn{x_a^K - w^*}
	+\frac{99n\chi^2M^2}{r}
	\right)
	+\frac{72n M^2}{r KT}
	\\&\aeq{uses the definition of $r_{yz}$ in \cref{eq:tau_x_eta_yz}}
	\frac{2}{K^2}
	\left(
	\frac{7nM^2}{r}
	+6r\sqn{x_a^K - w^*}
	+\frac{99n\chi^2M^2}{r}
	\right)
	+\frac{72n M^2}{r KT}
	\\&\leq
	\frac{212 n\chi^2 M^2}{rK^2}
	+\frac{72n M^2}{r KT}
	+\frac{12 r}{K^2}\sqn{x_a^K - w^*}
	\\&\aleq{uses the previously obtained upper-bound on $r\sqn{x_a^K - w^*}$}
	\frac{212 n\chi^2 M^2}{rK^2}
	+\frac{72n M^2}{r KT}
	+\frac{12}{K^2}\left(
	\frac{1386 n \chi^2 M^2}{rK^2}
	+\frac{72 n M^2}{r KT}
	\right),
\end{align*}
where \annotate.
Dividing both sides of the inequality by $n$ gives the following:
\begin{align*}
	p(x_o^K) - p(x^*)
	\leq
	\frac{212 \chi^2 M^2}{rK^2}
	+\frac{72 M^2}{r KT}
	+\frac{12}{K^2}\left(
	\frac{1386  \chi^2 M^2}{rK^2}
	+\frac{72  M^2}{r KT}
	\right).
\end{align*}
Hence, choosing the parameters $K$ and $T$ such that
\begin{align*}
	K \geq \cO\left( \frac{\chi M}{\sqrt{r\epsilon}}\right)
	\quad\text{and}\quad
	K \times T \geq \cO\left(\frac{M^2}{r\epsilon}\right)
\end{align*}
implies $p(x_o^K) - p(x^*) \leq \epsilon$, which concludes the proof.\qed

\newpage
\subsection{Proof of Theorem~\ref{thm:alg}}\label{sec:proof:alg}

With $r = 0$, the original problem~\eqref{eq:main} turns into the following problem:
\begin{equation}\label{eq:main_cvx}
	\min_{x\in \sX} \left[\bar{f}(x) = \frac{1}{n}\sum_{i=1}^{n} f_i(x)\right].
\end{equation}
Let $x^* \in \sX$ be the solution to problem~\eqref{eq:main_cvx}, such that $\norm{x^*} \leq R$, which always exists due to \Cref{ass:sol}. Let $r > 0$ be an arbitrary regularization parameter. We can upper-bound function $\bar{f}(x)$ using the regularized objective function $p(x)$ defined in \cref{eq:main} as follows:
\begin{align*}
	\bar{f}(x) \leq \bar{f}(x) + \frac{r}{2}\sqn{x} = p(x).
\end{align*}
On the other hand, we can lower-bound $\bar{f} (x^*)$ as follows:
\begin{align*}
	\bar{f}(x^*) = p(x^*) - \frac{r}{2}\sqn{x^*}
	\geq \min_{x' \in \sX} p(x') - \frac{r}{2}\sqn{x^*}
	\geq\min_{x' \in \sX} p(x') - \frac{rR^2}{2}.
\end{align*}
Hence, we can upper-bound the function suboptimality gap in problem~\eqref{eq:main_cvx} as follows:
\begin{align*}
	\bar{f}(x) - \bar{f}(x^*) \leq p(x) - \min_{x' \in \sX} p(x') + \frac{rR^2}{2}.
\end{align*}
Let the regularization parameter $r>0$ be chosen as follows:
\begin{equation}\label{eq:reg}
	r = {\epsilon}/{R^2}.
\end{equation}
Then, we obtain the following:
\begin{equation}\label{eq:aux2}
	\bar{f}(x) - \bar{f}(x^*) \leq p(x) - \min_{x' \in \sX} p(x') + \frac{\epsilon}{2}.
\end{equation}
We can apply \Cref{alg} to solving the regularized problem~\eqref{eq:main} with the regularization parameter $r$ defined in \cref{eq:reg}. \Cref{thm:alg_sc} implies that, to reach precision
\begin{equation}\label{eq:aux3}
	p(x_o^K) - \min_{x' \in \sX} p(x') \leq \frac{\epsilon}{2}
\end{equation}
it is sufficient to perform the following number of decentralized communications:
\begin{equation}
	K = \cO\left( \frac{\chi M}{\sqrt{r\epsilon}}\right)
	\aeq{use the definition of $r$ in \cref{eq:reg}}
	\cO\left( \frac{\chi MR}{\epsilon}\right),
\end{equation}
and the following number of subgradient computations:
\begin{equation}
	K\times T = \cO\left(\frac{M^2}{r\epsilon}\right)
	\aeq{use the definition of $r$ in \cref{eq:reg}}
	\cO\left(\frac{M^2R^2}{\epsilon^2}\right),
\end{equation}
where \annotate. Using \cref{eq:aux2,eq:aux3}, we also obtain the desired precision $\bar{f}(x_o^K) - \bar{f}(x^*) \leq \epsilon$, which concludes the proof. \qed

\newpage

\section{Proofs of Lemmas from Section~\ref{sec:alg_lem}}\label{sec:alg_lem_proof}

\subsection{Proof of Lemma~\ref{lem:sol}}\label{sec:proof:sol}

First, we pick the solution $x^* \in \sX$ to problem~\eqref{eq:main}, which is unique due to \Cref{ass:sol} and the fact that $r > 0$. Next, we define $w^* \in \cL$ as follows:
\begin{equation}\label{eq:w_star}
	w^* = (x^*,\ldots,x^*).
\end{equation}

From \Cref{ass:cvx,ass:lip} it follows that $\dom p(x) = \sX$ and $\dom f_i(x) = \sX$ for all $i \in \rng{1}{n}$, which implies the following:
\begin{equation}
	0 \in \partial p(x^*) =  r x^* + \frac{1}{n}\sum_{i=1}^{n} \partial f_i(x^*).
\end{equation}
Hence, there exists a vector $\Delta^* = (\Delta_1^*,\ldots,\Delta_n^*) \in (\sX)^n$ such that $\Delta_i^* \in \partial f_i(x^*)$ for all $i \in \rng{1}{n}$, and the following relation holds:
\begin{equation}\label{eq:aux1}
	rx^* + \frac{1}{n}\sum_{i=1}^{n} \Delta_i^* = 0.
\end{equation}

Next, we define $y^*\in (\sX)^n$ as follows:
\begin{equation}\label{eq:y_star}
	y^* = \Delta^* + r_x w^*.
\end{equation}
From \Cref{ass:cvx,ass:lip} it follows that $\dom F(x) = (\sX)^n$, which implies $y^* \in \partial F(x^*)$ and $0 \in \partial (F(\cdot) - \<y^*,\cdot>)(x^*) = \partial_x Q(x^*,y^*,z^*)$.

Next, we define $z^* \in \cL^\perp$ as follows:
\begin{equation}\label{eq:z_star}
	z^* = -rw^* - \Delta^*.
\end{equation}
Note that the inclusion $z^* \in \cL^\perp$ is implied by \cref{eq:aux1}.
Further, we get
\begin{align*}
	\nabla_z Q(w^*,y^*,z^*)
	\aeq{uses the definition of $Q(x,y,z)$ in \cref{eq:spp}}
	-r_{yz}(y^* + z^*)
	\aeq{uses the definition of $y^*$ and $z^*$}
	-r_{yz}(r_x-r)w^*
	\in \cL,
\end{align*}
where \annotate, and the last inclusion follows from the definition of $w^*$.
Moreover, we obtain the following
\begin{align*}
	\nabla_y Q(w^*,y^*,z^*)
	\aeq{uses the definition of $Q(x,y,z)$ in \cref{eq:spp}}
	-w^* - r_{yz}(y^* + z^*)
	\aeq{uses the definition of $y^*$ and $z^*$}
	-r_{yz}(r_{yz}^{-1} + r_x - r) w^*
	\aeq{uses \cref{eq:r_xyz}}
	0,
\end{align*}
where \annotate.

From \Cref{ass:lip} it follows that $\norm{\Delta_i^*} \leq M$ for all $i \in \rng{1}{n}$. Hence, using \cref{eq:aux1}, we get $r\norm{x^*} \leq M$, which implies $\sqn{w^*} \leq nM^2/r^2$. Moreover, we get
\begin{align*}
	\norm{y^*} \leq \norm{\Delta^*} + r_x\norm{w^*} \leq \sqrt{n}(M + r_xM/r) = (1+r_x / r)\sqrt{n}M,
\end{align*}
which implies $\sqn{y^*} \leq (1+r_x/r)^2 nM^2$.
Finally, we obtain
\begin{align*}
	\norm{z^*} \leq r\norm{w^*} + \norm{\Delta^*} \leq 2\sqrt{n}M,
\end{align*}
which implies $\sqn{z^*} \leq 4nM^2$ and concludes the proof.\qed

\newpage

\subsection{Proof of Lemma~\ref{lem:inner}}\label{sec:proof:inner}
We start with the following upper-bound on $\frac{1}{2\eta_x^k}\sqn{x^{k,t+1} - x}$:
\begin{align*}
	\mind{3em}\frac{1}{2\eta_x^k}\sqn{x^{k,t+1} - x}
	\\& \aeq{uses the parallelogram rule}
	\frac{1}{2\eta_x^k}\sqn{x^{k,t} - x} - \frac{1}{2\eta_x^k}\sqn{x^{k,t+1} - x^{k,t}}
	+\frac{1}{\eta_x^k}\<x^{k,t+1} - x^{k,t},x^{k,t+1} - x>
	\\&\aeq{uses \cref{line:ix} of \Cref{alg}}
	\frac{1}{2\eta_x^k}\sqn{x^{k,t} - x} - \frac{1}{2\eta_x^k}\sqn{x^{k,t+1} - x^{k,t}}
	\\&-\<g_x^{k,t} + \beta_k x^{k,t+1} - y^{k+1} + \tau_x^k(x^{k,t+1} - x^k),x^{k,t+1} - x>
	\\&=
	\frac{1}{2\eta_x^k}\sqn{x^{k,t} - x} - \frac{1}{2\eta_x^k}\sqn{x^{k,t+1} - x^{k,t}} + \<y^{k+1},x^{k,t+1} - x>
	\\&-\<\beta_k x^{k,t+1} + \tau_x^k(x^{k,t+1} - x^k),x^{k,t+1} - x>
	-\<g_x^{k,t},x^{k,t+1} - x>
	\\&\aleq{uses the parallelogram rule}
	\frac{1}{2\eta_x^k}\sqn{x^{k,t} - x} - \frac{1}{2\eta_x^k}\sqn{x^{k,t+1} - x^{k,t}} + \<y^{k+1},x^{k,t+1} - x>
	\\&
	-\frac{\tau_x^k}{2}\sqn{x^{k,t+1} - x^k}
	-\frac{\tau_x^k}{2}\sqn{x^{k,t+1} - x}
	+\frac{\tau_x^k}{2}\sqn{x^k - x}
	\\&
	-\frac{\beta_k}{2}\sqn{x^{k,t+1}}
	-\frac{\beta_k}{2}\sqn{x^{k,t+1} - x}
	+\frac{\beta_k}{2}\sqn{x}
	-\<g_x^{k,t},x^{k,t} - x> -\<g_x^{k,t},x^{k,t+1} - x^{k,t}>
	\\&\aleq{uses \cref{line:ig} of \Cref{alg}, \Cref{def:sub} and \Cref{ass:cvx}}
	\frac{1}{2\eta_x^k}\sqn{x^{k,t} - x} - \frac{1}{2\eta_x^k}\sqn{x^{k,t+1} - x^{k,t}} + \<y^{k+1},x^{k,t+1} - x>
	\\&
	-\frac{\tau_x^k}{2}\sqn{x^{k,t+1} - x^k}
	-\frac{\tau_x^k}{2}\sqn{x^{k,t+1} - x}
	+\frac{\tau_x^k}{2}\sqn{x^k - x}
	\\&
	-\frac{\beta_k}{2}\sqn{x^{k,t+1}}
	-\frac{\beta_k}{2}\sqn{x^{k,t+1} - x}
	+\frac{\beta_k}{2}\sqn{x}
	\\&
	+\sum_{i=1}^n(f_i(x_i) - f_i(x_i^{k,t})-\<g_{x,i}^{k,t},x_i^{k,t+1} - x_i^{k,t}>),
\end{align*}
where $(g_{x,1}^{k,t},\ldots,g_{x,n}^{k,t}) = (\subgrad f_1(x_1^{k,t}),\ldots,\subgrad f_n(x_n^{k,t})) =  g_x^{k,t} \in (\sX)^n$, \annotate.
Further, we obtain
\begin{align*}
	\mind{3em}\frac{1}{2\eta_x^k}\sqn{x^{k,t+1} - x}
	\\&\aleq{uses \Cref{ass:lip} and the Cauchy-Schwarz inequality}
	\frac{1}{2\eta_x^k}\sqn{x^{k,t} - x} - \frac{1}{2\eta_x^k}\sqn{x^{k,t+1} - x^{k,t}} + \<y^{k+1},x^{k,t+1} - x>
	\\&
	-\frac{\tau_x^k}{2}\sqn{x^{k,t+1} - x^k}
	-\frac{\tau_x^k}{2}\sqn{x^{k,t+1} - x}
	+\frac{\tau_x^k}{2}\sqn{x^k - x}
	\\&
	-\frac{\beta_k}{2}\sqn{x^{k,t+1}}
	-\frac{\beta_k}{2}\sqn{x^{k,t+1} - x}
	+\frac{\beta_k}{2}\sqn{x}
	\\&
	+\sum_{i=1}^n(f_i(x_i) - f_i(x_i^{k,t+1}) + M\norm{x_i^{k,t+1} - x_i^{k,t}}
	+\norm{g_{x,i}^{k,t}}\norm{x_i^{k,t+1} - x_i^{k,t}})
	\\&\aleq{uses the inequality $\norm{g_{x,i}^{k,t}} \leq M$, which follows from \Cref{ass:lip}}
	\frac{1}{2\eta_x^k}\sqn{x^{k,t} - x} - \frac{1}{2\eta_x^k}\sqn{x^{k,t+1} - x^{k,t}} + \<y^{k+1},x^{k,t+1} - x>
	\\&
	-\frac{\tau_x^k}{2}\sqn{x^{k,t+1} - x^k}
	-\frac{\tau_x^k}{2}\sqn{x^{k,t+1} - x}
	+\frac{\tau_x^k}{2}\sqn{x^k - x}
	\\&
	-\frac{\beta_k}{2}\sqn{x^{k,t+1}}
	-\frac{\beta_k}{2}\sqn{x^{k,t+1} - x}
	+\frac{\beta_k}{2}\sqn{x}
	\\&
	+\sum_{i=1}^n(f_i(x_i) - f_i(x_i^{k,t+1}) + 2M\norm{x_i^{k,t+1} - x_i^{k,t}})
	\\&\aleq{uses the definition of $\beta_k$ in \cref{eq:eta_x_beta_sigma_k}, the definition of $F(x)$ in \cref{eq:FG} and Young's inequality}
	\frac{1}{2\eta_x^k}\sqn{x^{k,t} - x} - \frac{1}{2\eta_x^k}\sqn{x^{k,t+1} - x^{k,t}} + \<y^{k+1},x^{k,t+1} - x>
	\\&
	-\frac{\tau_x^k}{2}\sqn{x^{k,t+1} - x^k}
	-\frac{\tau_x^k}{2}\sqn{x^{k,t+1} - x}
	+\frac{\tau_x^k}{2}\sqn{x^k - x}
	\\&
	+F(x)-F(x^{k,t+1})
	-\frac{r_x}{2}\sqn{x^{k,t+1} - x}
	+\sum_{i=1}^n\left(\frac{1}{2\eta_x^k}\sqn{x_i^{k,t+1} - x_i^{k,t}} + 2\eta_x^k M^2\right)
	\\&=
	\frac{1}{2\eta_x^k}\sqn{x^{k,t} - x}
	+\frac{\tau_x^k}{2}\sqn{x^k - x}
	+\<y^{k+1},x^{k,t+1} - x>
	+2n\eta_x^k M^2
	\\&
	-\frac{\tau_x^k}{2}\sqn{x^{k,t+1} - x^k}
	-\frac{\tau_x^k + r_x}{2}\sqn{x^{k,t+1} - x}
	-F(x^{k,t+1})
	+F(x),
\end{align*}
where \annotate.
After rearranging, we obtain
\begin{align*}
	\mind{3em}
	\frac{1}{2\eta_x^k}\sqn{x^{k,t+1} - x}
	\leq
	\frac{1}{2\eta_x^k}\sqn{x^{k,t} - x}
	+\frac{\tau_x^k}{2}\sqn{x^k - x}
	+2n\eta_x^k M^2
	-\Delta^{k,t+1},
\end{align*}
where $\Delta^{k,t+1}$ is defined as
\begin{align*}
	\Delta^{k,t+1}
	 & =
	F(x^{k,t+1}) - F(x)
	-\<y^{k+1},x^{k,t+1} - x>
	\\&
	+\frac{\tau_x^k + r_x}{2}\sqn{x^{k,t+1} - x}
	+\frac{\tau_x^k}{2}\sqn{x^{k,t+1} - x^k}
\end{align*}
Now, we sum these inequalities for $t=0,\ldots,T-1$ and obtain
\begin{align*}
	\frac{1}{2\eta_x^k}\sqn{x^{k,T} - x}
	\leq
	\frac{1}{2\eta_x^k}\sqn{x^{k,0} - x}
	+\frac{\tau_x^kT}{2}\sqn{x^k - x}
	+2n\eta_x^k M^2T
	-\sum_{t=1}^{T} \Delta^{k,t}.
\end{align*}
Dividing both sides of the inequality by $T$ gives
\begin{align*}
	\frac{1}{2\eta_x^kT}\sqn{x^{k,T} - x}
	\leq
	\frac{1}{2\eta_x^kT}\sqn{x^{k,0} - x}
	+\frac{\tau_x^k}{2}\sqn{x^k - x}
	+2n\eta_x^k M^2
	-\frac{1}{T}\sum_{t=1}^{T} \Delta^{k,t}.
\end{align*}
Using the definition of $\Delta^{k,t}$, the definition of $\tilde{x}^k$ on \cref{line:ox} of \Cref{alg} and \Cref{ass:cvx}, we obtain
\begin{align*}
	\frac{1}{2\eta_x^kT}\sqn{x^{k,T} - x}
	 & \leq
	\frac{1}{2\eta_x^kT}\sqn{x^{k,0} - x}
	+\frac{\tau_x^k}{2}\sqn{x^k - x}
	+2n\eta_x^k M^2
	\\&-
	\left(
	F(\tilde{x}^{k+1}) - F(x)
	-\<y^{k+1},\tilde{x}^{k+1} - x>
	\right)
	\\&
	-\left(
	\frac{\tau_x^k+r_x}{2}\sqn{\tilde{x}^{k+1} - x}
	+\frac{\tau_x^k}{2}\sqn{\tilde{x}^{k+1} - x^k}
	\right).
\end{align*}
Using the definition of $\eta_x^k$ and $\beta_k$ in \cref{eq:eta_x_beta_sigma_k}, we obtain
\begin{align*}
	\mind{10em}
	\frac{\tau_x^k}{2}\sqn{x^{k,T} - x}
	+\frac{\tau_x^k+\beta_k}{2}\sqn{\tilde{x}^{k+1} - x}
	\leq
	\frac{\tau_x^k}{2}\sqn{x^{k,0} - x}
	+\frac{\tau_x^k}{2}\sqn{x^k - x}
	+2n\eta_x^k M^2
	\\&-
	\left(
	F(\tilde{x}^{k+1}) - F(x)
	-\<y^{k+1},\tilde{x}^{k+1} - x>
	+\frac{\tau_x^k}{2}\sqn{\tilde{x}^{k+1} - x^k}
	\right).
\end{align*}
Using the definition of $x^{k,0}$ on \cref{line:iinit} of \Cref{alg}, the definition of $x^{k+1}$ on \cref{line:ox} of \Cref{alg}, the definition of $\eta_x^k$, $\beta_k$ and $\sigma_k$ in \cref{eq:eta_x_beta_sigma_k} and the convexity of $\norm{\cdot}$, we obtain
\begin{align*}
	(\tau_x^k + \tfrac{1}{2}r_x)\sqn{x^{k+1} - x}
	 & \leq
	\tau_x^k\sqn{x^k - x}
	+\frac{2n M^2}{\tau_x^kT}
	\\&-
	\left(
	F(\tilde{x}^{k+1}) - F(x)
	-\<y^{k+1},\tilde{x}^{k+1} - x>
	+\frac{\tau_x^k}{2}\sqn{\tilde{x}^{k+1} - x^k}
	\right),
\end{align*}
which concludes the proof. \qed

\newpage

\subsection{Proof of Lemma~\ref{lem:P}}\label{sec:proof:P}
Using \Cref{ass:gossip}, and the definition of $\mP$ in \cref{eq:P}, we obtain
\begin{equation}\label{eq:PW}
	\mP(\mW_k \otimes \mI_d) = (\mW_k \otimes \mI_d)\mP = (\mW_k \otimes \mI_d).
\end{equation}
Then, the desired relations can be trivially obtained by analyzing the lines of \Cref{alg}.\qed

\subsection{Proof of Lemma~\ref{lem:EF}}\label{sec:proof:EF}
We can upper-bound $\sqn{\eta_z^{k+1}m^{k+1}}_{\mP}$ as follows:
\begin{align*}
	\sqn{\eta_z^{k+1}m^{k+1}}_{\mP}
	 & \aeq{uses \cref{line:ext} of \Cref{alg}}
	\sqn{\eta_z^k(m^k + g_z^k - \hat{g}_z^k)}_{\mP}
	\\&\aeq{uses \cref{line:grad_comm} of \Cref{alg}}
	\sqn{\eta_z^k(m^k + g_z^k - (\mW_k\otimes \mI_d)(m^k + g_z^k))}_{\mP}
	\\&\aeq{uses \cref{eq:PW}}
	\sqn{\eta_z^k(\mP(m^k + g_z^k) - (\mW_k\otimes \mI_d)\mP(m^k + g_z^k))}
	\\&\aleq{uses \Cref{ass:chi}}
	(1 - 1/\chi)\sqn{\eta_z^k\mP(m^k + g_z^k)}
	\\&\aleq{uses the parallelogram rule and Young's inequality}
	(1 - 1/\chi)\left((1 + 1/(2\chi))\sqn{\eta_z^k m^k}_{\mP} + (1 + 2\chi)\sqn{\eta_z^k g_z^k}_{\mP}\right)
	\\&\leq
	(1-1/(2\chi))\sqn{\eta_z^k m^k}_{\mP}
	+2\chi\sqn{\eta_z^k g_z^k}_{\mP}
\end{align*}
where \annotate. Using this, we obtain
\begin{align*}
	\sqn{\eta_z^k m^k}_{\mP} \leq 2\chi\sqn{\eta_z^k m^k}_{\mP}
	-2\chi\sqn{\eta_z^{k+1} m^{k+1}}_{\mP}
	+4\chi^2\sqn{\eta_z^k g_z^k}_{\mP},
\end{align*}
which concludes the proof.\qed

\subsection{Proof of Lemma~\ref{lem:descent}}\label{sec:proof:descent}

We  can upper bound $\sqn{\tilde{g}_z^k - \mP g_z^k}$ as follows:
\begin{align*}
	\sqn{\tilde{g}_z^k - \mP g_z^k}
	 & \aeq{uses \cref{line:grad_comm} of \Cref{alg}}
	\sqn{(\mW_k \otimes \mI_d)g_z^k - \mP g_z^k}
	\\&\aeq{uses \cref{eq:PW}}
	\sqn{(\mW_k \otimes \mI_d)\mP g_z^k - \mP g_z^k}
	\\&\aleq{uses \Cref{ass:chi}}
	(1-1/\chi)\sqn{g_z^k}_{\mP}
\end{align*}
where \annotate.
On the other hand, $\sqn{\tilde{g}_z^k - \mP g_z^k}$ is equal to the following:
\begin{align*}
	\sqn{\tilde{g}_z^k - \mP g_z^k}
	 & \aeq{uses the parallelogram rule}
	\sqn{\tilde{g}_z^k}
	+\sqn{g_z^k}_\mP
	-2\<\tilde{g}_z^k,\mP g_z^k>
	\\&\aeq{uses \cref{line:ext} of \Cref{alg}}
	\frac{1}{(\theta_z^k)^2}\sqn{\ol{z}^{k+1} - \uln{z}^k}
	+\sqn{g_z^k}_\mP
	+\frac{2}{\theta_z^k}\<\ol{z}^{k+1} - \uln{z}^k,\mP g_z^k>.
\end{align*}
where \annotate.
Hence, we obtain the following
\begin{align*}
	\frac{1}{(\theta_z^k)^2}\sqn{\ol{z}^{k+1} - \uln{z}^k}
	+\frac{2}{\theta_z^k}\<\ol{z}^{k+1} - \uln{z}^k,\mP g_z^k>
	+\frac{1}{\chi}\sqn{g_z^k}_{\mP}
	\leq
	0.
\end{align*}
After rearranging and multiplying both sides of the inequality by $\frac{\theta_z^k}{2\alpha_k}$, we obtain
\begin{align*}
	0 & \geq
	\alpha_k^{-1}\left(
	\<\ol{z}^{k+1} - \uln{z}^k,\mP g_z^k>
	+\frac{1}{2\theta_z^k}\sqn{\ol{z}^{k+1} - \uln{z}^k}
	\right)
	+\frac{\theta_z^k}{2\alpha_k\chi}\sqn{g_z^k}_{\mP}
	\\&\aeq{uses \cref{eq:theta_z_k}}
	\alpha_k^{-1}\left(
	\<\mP(\ol{z}^{k+1} - \uln{z}^k),g_z^k>
	+r_{yz}\sqn{\ol{z}^{k+1} - \uln{z}^k}
	\right)
	+\frac{1}{4\alpha_k\chi r_{yz}}\sqn{g_z^k}_{\mP}
	\\&\aeq{uses \Cref{lem:P}}
	\alpha_k^{-1}\left(
	\<\ol{z}^{k+1} - \uln{z}^k,g_z^k>
	+r_{yz}\sqn{\ol{z}^{k+1} - \uln{z}^k}
	\right)
	+\frac{1}{4\alpha_k\chi r_{yz}}\sqn{g_z^k}_{\mP}
\end{align*}
where \annotate, which concludes the proof. \qed

\newpage
\subsection{Proof of Lemma~\ref{lem:outer}}\label{sec:proof:outer}

We can upper-bound $\frac{1}{2\eta_y^k}\sqn{y^{k+1} - y}$ as follows:
\begin{align*}
	\frac{1}{2\eta_y^k}\sqn{y^{k+1} - y}
	 & \aeq{uses the parallelogram rule}
	\frac{1}{2\eta_y^k}\sqn{y^{k} - y}
	-\frac{1}{2\eta_y^k}\sqn{y^{k+1} - y^k}
	+\frac{1}{\eta_y^k}\<y^{k+1} - y^k, y^{k+1} - y>
	\\&\aeq{uses \cref{line:yz} of \Cref{alg}}
	\frac{1}{2\eta_y^k}\sqn{y^{k+1} - y}
	-\frac{1}{2\eta_y^k}\sqn{y^{k+1} - y^k}
	-\<g_y^k + \hat{x}^{k+1}, y^{k+1} - y>
	\\&=
	\frac{1}{2\eta_y^k}\sqn{y^{k} - y}
	-\frac{1}{2\eta_y^k}\sqn{y^{k+1} - y^k}
	-\<\hat{x}^{k+1}, y^{k+1} - y>
	\\&
	-\<g_y^k, y^{k+1} - y^k + y^k - \uln{y}^k + \uln{y}^k - y>
	\\&\aeq{uses \Cref{line:comb,line:ext} of \Cref{alg}}
	\frac{1}{2\eta_y^k}\sqn{y^{k} - y}
	-\frac{1}{2\eta_y^k}\sqn{y^{k+1} - y^k}
	-\<\hat{x}^{k+1}, y^{k+1} - y>
	\\&
	-\alpha_k^{-1}\<g_y^k, \ol{y}^{k+1} - \uln{y}^k>
	+(1-\alpha_k)\alpha_k^{-1}\<g_y^k,\ol{y}^k - \uln{y}^k>
	+\<g_y^k,y-\uln{y}^k>
	\\&\aeq{uses \cref{line:yz} of \Cref{alg}}
	\frac{1}{2\eta_y^k}\sqn{y^{k} - y}
	-\frac{1}{2\eta_y^k}\sqn{y^{k+1} - y^k}
	-\<\tilde{x}^{k+1}, y^{k+1} - y>
	\\&
	-\alpha_k^{-1}\<g_y^k, \ol{y}^{k+1} - \uln{y}^k>
	+(1-\alpha_k)\alpha_k^{-1}\<g_y^k,\ol{y}^k - \uln{y}^k>
	+\<g_y^k,y-\uln{y}^k>
	\\&
	+\<\tilde{x}^{k+1} - \hat{x}^{k+1}, y^{k+1} - y>
\end{align*}
where \annotate.
Further, we can upper-bound the term $\<\tilde{x}^{k+1} - \hat{x}^{k+1}, y^{k+1} - y>$ as follows:
\begin{align*}
	\mind{3em}
	\<\tilde{x}^{k+1} - \hat{x}^{k+1}, y^{k+1} - y>
	\\& \aeq{uses \Cref{line:yz} of \Cref{alg}}
	\<\tilde{x}^{k+1} - x^k - \gamma_k (\tilde{x}^k - x^{k-1}), y^{k+1} - y>
	\\&=
	\gamma_k\<x^{k-1} - \tilde{x}^k, y^k - y>
	-\<x^k - \tilde{x}^{k+1}, y^{k+1} - y>
	+\gamma_k\<x^{k-1} - \tilde{x}^k, y^{k+1} - y^k>
	\\&\aleq{uses Young's inequality}
	\gamma_k\<x^{k-1} - \tilde{x}^k, y^k - y>
	-\<x^k - \tilde{x}^{k+1}, y^{k+1} - y>
	+\frac{1}{4\eta_y^k}\sqn{y^{k+1}-y^k}
	\\&
	+2\eta_y^{k}\gamma_k^2\sqn{x^{k-1} - \tilde{x}^k}.
\end{align*}
where \annotate.
Plugging this into the previous inequality gives
\begin{align*}
	\frac{1}{2\eta_y^k}\sqn{y^{k+1} - y}
	 & \leq
	\frac{1}{2\eta_y^k}\sqn{y^{k} - y}
	-\frac{1}{4\eta_y^k}\sqn{y^{k+1} - y^k}
	-\<\tilde{x}^{k+1}, y^{k+1} - y>
	\\&
	-\alpha_k^{-1}\<g_y^k, \ol{y}^{k+1} - \uln{y}^k>
	+(1-\alpha_k)\alpha_k^{-1}\<g_y^k,\ol{y}^k - \uln{y}^k>
	+\<g_y^k,y-\uln{y}^k>
	\\&
	+\gamma_k\<x^{k-1} - \tilde{x}^k, y^k - y>
	-\<x^k - \tilde{x}^{k+1}, y^{k+1} - y>
	+2\eta_y^{k}\gamma_k^2\sqn{x^{k-1} - \tilde{x}^k}
	\\&\aeq{\cref{line:ext} of \Cref{alg}}
	\frac{1}{2\eta_y^k}\sqn{y^{k} - y}
	-\<\tilde{x}^{k+1}, y^{k+1} - y>
	+2\eta_y^{k}\gamma_k^2\sqn{x^{k-1} - \tilde{x}^k}
	\\&
	+\<g_y^k,y-\uln{y}^k>
	+(1-\alpha_k)\alpha_k^{-1}\<g_y^k,\ol{y}^k - \uln{y}^k>
	-\alpha_k^{-1}\<g_y^k, \ol{y}^{k+1} - \uln{y}^k>
	\\&
	-\frac{1}{4\eta_y^k\alpha_k^2}\sqn{\ol{y}^{k+1} - \uln{y}^k}
	+\gamma_k\<x^{k-1} - \tilde{x}^k, y^k - y>
	-\<x^k - \tilde{x}^{k+1}, y^{k+1} - y>
	\\&\aeq{uses \cref{eq:tau_x_eta_yz_k,eq:tau_x_eta_yz}}
	\frac{1}{2\eta_y^k}\sqn{y^{k} - y}
	+2\eta_y^{k}\gamma_k^2\sqn{x^{k-1} - \tilde{x}^k}
	-\<\tilde{x}^{k+1}, y^{k+1} - y>
	\\&
	+\gamma_k\<x^{k-1} - \tilde{x}^k, y^k - y>
	-\<x^k - \tilde{x}^{k+1}, y^{k+1} - y>
	+\<g_y^k,y-\uln{y}^k>
	\\&
	+(1-\alpha_k)\alpha_k^{-1}\<g_y^k,\ol{y}^k - \uln{y}^k>
	-\alpha_k^{-1}\left(\<g_y^k, \ol{y}^{k+1} - \uln{y}^k> + r_{yz}\sqn{\ol{y}^{k+1} - \uln{y}^k}\right),
\end{align*}
where \annotate.

Let $\hat{z}^k$ be defined for all $k\in\rng{0}{K}$ as follows:
\begin{equation}\label{eq:z_hat}
	\hat{z}^k = z^k - \eta_z^k \mP m^k.
\end{equation}
Using \cref{eq:z_hat} and \cref{line:yz,line:ext} of \Cref{alg}, we obtain
\begin{align*}
	\hat{z}^{k+1} & = \hat{z}^k + z^{k+1} - z^k - \mP(\eta_z^{k+1}m^{k+1} - \eta_z^k m^k)
	\\&=
	\hat{z}^k
	-\mP(\eta_z^k\hat{g}_z^k + \eta_z^{k+1}m^{k+1} - \eta_z^k m^k)
	\\&=
	\hat{z}^k
	-\mP(\eta_z^k\hat{g}_z^k + \eta_z^{k+1}(\eta_z^k/\eta_z^{k+1})(m^k +g_z^k  - \hat{g}_z^k) - \eta_z^k m^k)
	\\&=
	\hat{z}^k-\eta_z^k\mP g_z^k.
\end{align*}
Hence, we can upper-bound $\frac{1}{2\eta_z^k}\sqn{\hat{z}^{k+1} - z}$ as follows:
\begin{align*}
	\frac{1}{2\eta_z^k}\sqn{\hat{z}^{k+1} - z}
	 & \aeq{uses the parallelogram rule}
	\frac{1}{2\eta_z^k}\sqn{\hat{z}^{k} - z}
	+\frac{1}{2\eta_z^k}\sqn{\hat{z}^{k+1} - \hat{z}^k}
	+\frac{1}{\eta_z^k}\<\hat{z}^{k+1} - \hat{z}^k,\hat{z}^{k} - z>
	\\&\aeq{uses the update rule for $\hat{z}^k$ which we previously obtained}
	\frac{1}{2\eta_z^k}\sqn{\hat{z}^{k} - z}
	+\frac{\eta_z^k}{2}\sqn{g_z^k}_{\mP}
	-\<\mP g_z^k,\hat{z}^{k} - z>
	\\&=
	\frac{1}{2\eta_z^k}\sqn{\hat{z}^{k} - z}
	+\frac{\eta_z^k}{2}\sqn{g_z^k}_{\mP}
	-\<\mP g_z^k,z^k - z>
	+\<\mP g_z^k,z^k - \hat{z}^{k}>
	\\&\aeq{uses \cref{eq:z_hat}}
	\frac{1}{2\eta_z^k}\sqn{\hat{z}^{k} - z}
	+\frac{\eta_z^k}{2}\sqn{g_z^k}_{\mP}
	-\<\mP g_z^k,z^k - z>
	+\eta_z^k\<\mP g_z^k,\mP m^k>,
\end{align*}
where \annotate.
Further, we can upper-bound the term $\eta_z^k\<\mP g_z^k,\mP m^k>$ as follows
\begin{align*}
	\eta_z^k\<\mP g_z^k,\mP m^k>
	 & \aleq{uses the Cauchy-Schwarz inequality}
	\frac{1}{\eta_z^k}\norm{\eta_z^kg_z^k}_{\mP}\norm{\eta_z^km^k}_{\mP}
	\\&\aleq{uses Young's inequality}
	\frac{1}{2\eta_z^k}\left(2\chi\sqn{\eta_z^kg_z^k}_{\mP} + \frac{1}{2\chi}\sqn{\eta_z^km^k}_{\mP}\right)
	\\&\aleq{uses \Cref{lem:EF}}
	\frac{1}{2\eta_z^k}\left(
	4\chi\sqn{\eta_z^kg_z^k}_{\mP}
	+
	\sqn{\eta_z^k m^k}_{\mP}
	-\sqn{\eta_z^{k+1} m^{k+1}}_{\mP}
	\right)
\end{align*}
where \annotate. Plugging this into the previous inequality gives
\begin{align*}
	\frac{1}{2\eta_z^k}\sqn{\hat{z}^{k+1} - z}
	 & \leq
	\frac{1}{2\eta_z^k}\sqn{\hat{z}^{k} - z}
	+\frac{\eta_z^k}{2}\sqn{g_z^k}_{\mP}
	-\<\mP g_z^k,z^k - z>
	\\&
	+\frac{1}{2\eta_z^k}\left(
	4\chi\sqn{\eta_z^kg_z^k}_{\mP}
	+
	\sqn{\eta_z^k m^k}_{\mP}
	-\sqn{\eta_z^{k+1} m^{k+1}}_{\mP}
	\right)
	\\&\aeq{uses \Cref{lem:P} and the fact that $z \in \cL^\perp$}
	\frac{1}{2\eta_z^k}\sqn{\hat{z}^{k} - z}
	+\frac{\eta_z^k(1+4\chi)}{2}\sqn{g_z^k}_{\mP}
	-\<g_z^k,z^k - \uln{z}^k + \uln{z}^k - z>
	\\&
	+\frac{1}{2\eta_z^k}\left(
	\sqn{\eta_z^k m^k}_{\mP}
	-\sqn{\eta_z^{k+1} m^{k+1}}_{\mP}
	\right)
	\\&\aeq{uses \cref{line:comb} of \Cref{alg}}
	\frac{1}{2\eta_z^k}\sqn{\hat{z}^{k} - z}
	+\frac{\eta_z^k(1+4\chi)}{2}\sqn{g_z^k}_{\mP}
	+\frac{1}{2\eta_z^k}\left(
	\sqn{\eta_z^k m^k}_{\mP}
	-\sqn{\eta_z^{k+1} m^{k+1}}_{\mP}
	\right)
	\\&
	+\<g_z^k,z - \uln{z}^k>
	+(1-\alpha_k)\alpha_k^{-1}\<g_z^k,\ol{z}^k - \uln{z}^k>
	\\&\aleq{uses \Cref{lem:descent}}
	\frac{1}{2\eta_z^k}\sqn{\hat{z}^{k} - z}
	+\frac{\eta_z^k(1+4\chi)}{2}\sqn{g_z^k}_{\mP}
	+\frac{1}{2\eta_z^k}\left(
	\sqn{\eta_z^k m^k}_{\mP}
	-\sqn{\eta_z^{k+1} m^{k+1}}_{\mP}
	\right)
	\\&
	+\<g_z^k,z - \uln{z}^k>
	+(1-\alpha_k)\alpha_k^{-1}\<g_z^k,\ol{z}^k - \uln{z}^k>
	\\&
	-\alpha_k^{-1}\left(
	\<\ol{z}^{k+1} - \uln{z}^k,g_z^k>
	+r_{yz}\sqn{\ol{z}^{k+1} - \uln{z}^k}
	\right)
	-\frac{1}{4\alpha_k\chi r_{yz}}\sqn{g_z^k}_{\mP}
	\\&\aleq{uses \cref{eq:tau_x_eta_yz_k,eq:tau_x_eta_yz}}
	\frac{1}{2\eta_z^k}\sqn{\hat{z}^{k} - z}
	+\frac{1}{2\eta_z^k}\left(
	\sqn{\eta_z^k m^k}_{\mP}
	-\sqn{\eta_z^{k+1} m^{k+1}}_{\mP}
	\right)
	+\<g_z^k,z - \uln{z}^k>
	\\&
	+(1-\alpha_k)\alpha_k^{-1}\<g_z^k,\ol{z}^k - \uln{z}^k>
	-\alpha_k^{-1}\left(
	\<\ol{z}^{k+1} - \uln{z}^k,g_z^k>
	+r_{yz}\sqn{\ol{z}^{k+1} - \uln{z}^k}
	\right)
\end{align*}
where \annotate.

Now we combine the upper-bounds for $\frac{1}{2\eta_y^k}\sqn{y^{k+1} - y}$ and $\frac{1}{2\eta_z^k}\sqn{\hat{z}^{k+1} - z}$ and obtain the following:
\begin{align*}
	\mind{3em}
	\frac{1}{2\eta_y^k}\sqn{y^{k+1} - y}
	+\frac{1}{2\eta_z^k}\sqn{\hat{z}^{k+1} - z}
	\\&\leq
	\frac{1}{2\eta_y^k}\sqn{y^{k} - y}
	+\frac{1}{2\eta_z^k}\sqn{\hat{z}^{k} - z}
	+2\eta_y^{k}\gamma_k^2\sqn{x^{k-1} - \tilde{x}^k}
	-\<\tilde{x}^{k+1}, y^{k+1} - y>
	\\&
	+\gamma_k\<x^{k-1} - \tilde{x}^k, y^k - y>
	-\<x^k - \tilde{x}^{k+1}, y^{k+1} - y>
	+\frac{1}{2\eta_z^k}\left(
	\sqn{\eta_z^k m^k}_{\mP}
	-\sqn{\eta_z^{k+1} m^{k+1}}_{\mP}
	\right)
	\\&
	+\<g_y^k,y-\uln{y}^k>+\<g_z^k,z - \uln{z}^k>
	+(1-\alpha_k)\alpha_k^{-1}\left(\<g_y^k,\ol{y}^k - \uln{y}^k> + \<g_z^k,\ol{z}^k - \uln{z}^k>\right)
	\\&
	-\alpha_k^{-1}\left(
	\<g_y^k, \ol{y}^{k+1} - \uln{y}^k> + \<\ol{z}^{k+1} - \uln{z}^k,g_z^k>
	+ r_{yz}\sqn{\ol{y}^{k+1} - \uln{y}^k}+r_{yz}\sqn{\ol{z}^{k+1} - \uln{z}^k}\right)
	\\&\aleq{uses the definition of $g_y^k$ and $g_z^k$ on \cref{line:grad} of \Cref{alg} and the convexity and $(2r_{yz})$-smoothness of the function $G(y,z)$}
	\frac{1}{2\eta_y^k}\sqn{y^{k} - y}
	+\frac{1}{2\eta_z^k}\sqn{\hat{z}^{k} - z}
	+2\eta_y^{k}\gamma_k^2\sqn{x^{k-1} - \tilde{x}^k}
	-\<\tilde{x}^{k+1}, y^{k+1} - y>
	\\&
	+\gamma_k\<x^{k-1} - \tilde{x}^k, y^k - y>
	-\<x^k - \tilde{x}^{k+1}, y^{k+1} - y>
	+\frac{1}{2\eta_z^k}\left(
	\sqn{\eta_z^k m^k}_{\mP}
	-\sqn{\eta_z^{k+1} m^{k+1}}_{\mP}
	\right)
	\\&
	+G(y,z) - G(\uln{y}^k,\uln{z}^k)
	+(1-\alpha_k)\alpha_k^{-1}\left(G(\ol{y}^k,\ol{z}^k) - G(\uln{y}^k,\uln{z}^k)\right)
	\\&
	-\alpha_k^{-1}\left(
	G(\ol{y}^{k+1},\ol{z}^{k+1})-G(\uln{y}^k,\uln{z}^k)
	\right)
	\\&=
	\frac{1}{2\eta_y^k}\sqn{y^{k} - y}
	+\frac{1}{2\eta_z^k}\sqn{\hat{z}^{k} - z}
	+2\eta_y^{k}\gamma_k^2\sqn{x^{k-1} - \tilde{x}^k}
	-\<\tilde{x}^{k+1}, y^{k+1} - y>
	\\&
	+\gamma_k\<x^{k-1} - \tilde{x}^k, y^k - y>
	-\<x^k - \tilde{x}^{k+1}, y^{k+1} - y>
	+\frac{1}{2\eta_z^k}\left(
	\sqn{\eta_z^k m^k}_{\mP}
	-\sqn{\eta_z^{k+1} m^{k+1}}_{\mP}
	\right)
	\\&
	+(1-\alpha_k)\alpha_k^{-1}\left(G(\ol{y}^k,\ol{z}^k) - G(y,z)\right)
	-\alpha_k^{-1}\left(
	G(\ol{y}^{k+1},\ol{z}^{k+1})-G(y,z)
	\right),
\end{align*}
where \annotate.
Further, we divide both sides of the inequality by $\alpha_k$ and, using \cref{eq:tau_x_eta_yz_k}, obtain the following:
\begin{align*}
	\mind{3em}
	\frac{1}{2\eta_y}\sqn{y^{k+1} - y}
	+\frac{1}{2\eta_z}\sqn{\hat{z}^{k+1} - z}
	+\frac{1}{2\eta_z}\sqn{\eta_z^{k+1} m^{k+1}}_{\mP}
	\\&\leq
	\frac{1}{2\eta_y}\sqn{y^{k} - y}
	+\frac{1}{2\eta_z}\sqn{\hat{z}^{k} - z}
	+\frac{1}{2\eta_z}\sqn{\eta_z^k m^k}_{\mP}
	+2\eta_y\alpha_k^{-2}\gamma_k^2\sqn{x^{k-1} - \tilde{x}^k}
	\\&
	+\gamma_k\alpha_k^{-1}\<x^{k-1} - \tilde{x}^k, y^k - y>
	-\alpha_k^{-1}\<x^k - \tilde{x}^{k+1}, y^{k+1} - y>
	-\alpha_k^{-1}\<\tilde{x}^{k+1}, y^{k+1} - y>
	\\&
	+(\alpha_k^{-2}-\alpha_k^{-1})\left(G(\ol{y}^k,\ol{z}^k) - G(y,z)\right)
	-\alpha_k^{-2}\left(
	G(\ol{y}^{k+1},\ol{z}^{k+1})-G(y,z)
	\right)
\end{align*}

Next, we divide the inequality in \Cref{lem:inner} by $\alpha_k$ and, using the definition of $\tau_x^k$ and $\tau_x$ in \cref{eq:tau_x_eta_yz_k,eq:tau_x_eta_yz}, obtain the following:
\begin{align*}
	\tau_x(\alpha_k^{-2} + \alpha_k^{-1})\sqn{x^{k+1} - x}
	 & \leq
	\tau_x\alpha_k^{-2}\sqn{x^k - x}
	-\frac{\tau_x\alpha_k^{-2}}{2}\sqn{\tilde{x}^{k+1} - x^k}
	+\frac{2n M^2}{\tau_xT}
	\\&-
	\alpha_k^{-1}\left(
	F(\tilde{x}^{k+1}) - F(x)
	-\<y^{k+1},\tilde{x}^{k+1} - x>
	\right).
\end{align*}
Combining this inequality with the previous upper-bound gives the following:
\begin{align*}
	\mind{2em}
	\tau_x(\alpha_k^{-2} + \alpha_k^{-1})\sqn{x^{k+1} - x}
	+\frac{1}{2\eta_y}\sqn{y^{k+1} - y}
	+\frac{1}{2\eta_z}\sqn{\hat{z}^{k+1} - z}
	+\frac{1}{2\eta_z}\sqn{\eta_z^{k+1} m^{k+1}}_{\mP}
	\\&\leq
	\tau_x\alpha_k^{-2}\sqn{x^k - x}
	+\frac{1}{2\eta_y}\sqn{y^{k} - y}
	+\frac{1}{2\eta_z}\sqn{\hat{z}^{k} - z}
	+\frac{1}{2\eta_z}\sqn{\eta_z^k m^k}_{\mP}
	+\frac{2n M^2}{\tau_xT}
	\\&
	-\frac{\tau_x\alpha_k^{-2}}{2}\sqn{\tilde{x}^{k+1} - x^k}
	+2\eta_y\alpha_k^{-2}\gamma_k^2\sqn{x^{k-1} - \tilde{x}^k}
	-\alpha_k^{-1}\<x^k - \tilde{x}^{k+1}, y^{k+1} - y>
	\\&
	+\gamma_k\alpha_k^{-1}\<x^{k-1} - \tilde{x}^k, y^k - y>
	-\alpha_k^{-1}\left(
	F(\tilde{x}^{k+1}) - F(x)
	+\<y^{k+1},x>
	-\<\tilde{x}^{k+1},y>
	\right)
	\\&
	+(\alpha_k^{-2}-\alpha_k^{-1})\left(G(\ol{y}^k,\ol{z}^k) - G(y,z)\right)
	-\alpha_k^{-2}\left(
	G(\ol{y}^{k+1},\ol{z}^{k+1})-G(y,z)
	\right)
	\\&\aeq{uses the fact that $y^{k+1} = \alpha_k^{-1}\ol{y}^{k+1} - (1-\alpha_k)\alpha_k^{-1}\ol{y}^k$, which follows from \cref{line:comb,line:ext} of \Cref{alg}}
	\tau_x\alpha_k^{-2}\sqn{x^k - x}
	+\frac{1}{2\eta_y}\sqn{y^{k} - y}
	+\frac{1}{2\eta_z}\sqn{\hat{z}^{k} - z}
	+\frac{1}{2\eta_z}\sqn{\eta_z^k m^k}_{\mP}
	+\frac{2n M^2}{\tau_xT}
	\\&
	-\frac{\tau_x\alpha_k^{-2}}{2}\sqn{\tilde{x}^{k+1} - x^k}
	+2\eta_y\alpha_k^{-2}\gamma_k^2\sqn{x^{k-1} - \tilde{x}^k}
	-\alpha_k^{-1}\<x^k - \tilde{x}^{k+1}, y^{k+1} - y>
	\\&
	+\gamma_k\alpha_k^{-1}\<x^{k-1} - \tilde{x}^k, y^k - y>
	-\alpha_k^{-1}\left(
	F(\tilde{x}^{k+1})
	-\<\tilde{x}^{k+1},y>
	-G(y,z)
	\right)
	\\&
	+(\alpha_k^{-2}-\alpha_k^{-1})\left(
	G(\ol{y}^k,\ol{z}^k)+\<\ol{y}^k,x> - F(x)
	\right)
	-\alpha_k^{-2}\left(
	G(\ol{y}^{k+1},\ol{z}^{k+1})+\<\ol{y}^{k+1},x> - F(x)
	\right)
	\\&\aeq{uses the definition of $Q(x,y,z)$ in \cref{eq:spp}}
	\tau_x\alpha_k^{-2}\sqn{x^k - x}
	+\frac{1}{2\eta_y}\sqn{y^{k} - y}
	+\frac{1}{2\eta_z}\sqn{\hat{z}^{k} - z}
	+\frac{1}{2\eta_z}\sqn{\eta_z^k m^k}_{\mP}
	+\frac{2n M^2}{\tau_xT}
	\\&
	-\frac{\tau_x\alpha_k^{-2}}{2}\sqn{\tilde{x}^{k+1} - x^k}
	+2\eta_y\alpha_k^{-2}\gamma_k^2\sqn{x^{k-1} - \tilde{x}^k}
	-\alpha_k^{-1}\<x^k - \tilde{x}^{k+1}, y^{k+1} - y>
	\\&
	+\gamma_k\alpha_k^{-1}\<x^{k-1} - \tilde{x}^k, y^k - y>
	-(\alpha_k^{-2}-\alpha_k^{-1})Q(x,\ol{y}^k,\ol{z}^k)
	+\alpha_k^{-2}Q(x,\ol{y}^{k+1},\ol{z}^{k+1})
	\\&
	-\alpha_k^{-1}Q(\tilde{x}^{k+1},y,z)
	\\&\aleq{uses \cref{line:ox} of \Cref{alg} and \Cref{ass:cvx}}
	\tau_x\alpha_k^{-2}\sqn{x^k - x}
	+\frac{1}{2\eta_y}\sqn{y^{k} - y}
	+\frac{1}{2\eta_z}\sqn{\hat{z}^{k} - z}
	+\frac{1}{2\eta_z}\sqn{\eta_z^k m^k}_{\mP}
	+\frac{2n M^2}{\tau_xT}
	\\&
	-\frac{\tau_x\alpha_k^{-2}}{2}\sqn{\tilde{x}^{k+1} - x^k}
	+2\eta_y\alpha_k^{-2}\gamma_k^2\sqn{x^{k-1} - \tilde{x}^k}
	-\alpha_k^{-1}\<x^k - \tilde{x}^{k+1}, y^{k+1} - y>
	\\&
	+\gamma_k\alpha_k^{-1}\<x^{k-1} - \tilde{x}^k, y^k - y>
	-(\alpha_k^{-2}-\alpha_k^{-1})Q(x,\ol{y}^k,\ol{z}^k)
	+\alpha_k^{-2}Q(x,\ol{y}^{k+1},\ol{z}^{k+1})
	\\&
	-\alpha_k^{-2}Q(\ol{x}^{k+1},y,z)
	+(\alpha_k^{-2} - \alpha_k^{-1})Q(\ol{x}^{k},y,z)
	\\&=
	\tau_x\alpha_k^{-2}\sqn{x^k - x}
	+\frac{1}{2\eta_y}\sqn{y^{k} - y}
	+\frac{1}{2\eta_z}\sqn{\hat{z}^{k} - z}
	+\frac{1}{2\eta_z}\sqn{\eta_z^k m^k}_{\mP}
	+\frac{2n M^2}{\tau_xT}
	\\&
	-\frac{\tau_x\alpha_k^{-2}}{2}\sqn{\tilde{x}^{k+1} - x^k}
	+2\eta_y\alpha_k^{-2}\gamma_k^2\sqn{x^{k-1} - \tilde{x}^k}
	-\alpha_k^{-1}\<x^k - \tilde{x}^{k+1}, y^{k+1} - y>
	\\&
	+\gamma_k\alpha_k^{-1}\<x^{k-1} - \tilde{x}^k, y^k - y>
	+(\alpha_k^{-2}-\alpha_k^{-1})\left(Q(\ol{x}^{k},y,z) - Q(x,\ol{y}^k,\ol{z}^k)\right)
	\\&
	-\alpha_k^{-2}\left(Q(\ol{x}^{k+1},y,z) - Q(x,\ol{y}^{k+1},\ol{z}^{k+1})\right)
\end{align*}
where \annotate.

Further, let $\alpha_K = 3 / (K+3)$. Then from \cref{eq:alpha_gamma_k} it follows that $\alpha_{k}^{-2} +\alpha_{k}^{-1} \geq \alpha_{k+1}^{-2}$, $\gamma_k \alpha_k^{-1} = (k+2)/3$ and $\alpha_k^{-1} = (k+3)/3$. Hence, we obtain the following:
\begin{align*}
	\mind{1.5em}
	\tau_x\alpha_{k+1}^{-2}\sqn{x^{k+1} - x}
	+\frac{1}{2\eta_y}\sqn{y^{k+1} - y}
	+\frac{1}{2\eta_z}\sqn{\hat{z}^{k+1} - z}
	+\frac{1}{2\eta_z}\sqn{\eta_z^{k+1} m^{k+1}}_{\mP}
	\\&\leq
	\tau_x\alpha_k^{-2}\sqn{x^k - x}
	+\frac{1}{2\eta_y}\sqn{y^{k} - y}
	+\frac{1}{2\eta_z}\sqn{\hat{z}^{k} - z}
	+\frac{1}{2\eta_z}\sqn{\eta_z^k m^k}_{\mP}
	+\frac{2n M^2}{\tau_xT}
	\\&
	-\frac{\tau_x(k+3)^2}{18}\sqn{\tilde{x}^{k+1} - x^k}
	+\frac{4\eta_y(k+2)^2}{18}\sqn{x^{k-1} - \tilde{x}^k}
	-\frac{k+3}{3}\<x^k - \tilde{x}^{k+1}, y^{k+1} - y>
	\\&
	+\frac{k+2}{3}\<x^{k-1} - \tilde{x}^k, y^k - y>
	+(\alpha_k^{-2}-\alpha_k^{-1})\left(Q(\ol{x}^{k},y,z) - Q(x,\ol{y}^k,\ol{z}^k)\right)
	\\&
	-\alpha_k^{-2}\left(Q(\ol{x}^{k+1},y,z) - Q(x,\ol{y}^{k+1},\ol{z}^{k+1})\right)
	\\&\aeq{uses \cref{eq:tau_x_eta_yz}}
	\tau_x\alpha_k^{-2}\sqn{x^k - x}
	+\frac{1}{2\eta_y}\sqn{y^{k} - y}
	+\frac{1}{2\eta_z}\sqn{\hat{z}^{k} - z}
	+\frac{1}{2\eta_z}\sqn{\eta_z^k m^k}_{\mP}
	+\frac{2n M^2}{\tau_xT}
	\\&
	-\frac{2\eta_y(k+3)^2}{9}\sqn{\tilde{x}^{k+1} - x^k}
	+\frac{2\eta_y(k+2)^2}{9}\sqn{x^{k-1} - \tilde{x}^k}
	-\frac{k+3}{3}\<x^k - \tilde{x}^{k+1}, y^{k+1} - y>
	\\&
	+\frac{k+2}{3}\<x^{k-1} - \tilde{x}^k, y^k - y>
	+(\alpha_k^{-2}-\alpha_k^{-1})\left(Q(\ol{x}^{k},y,z) - Q(x,\ol{y}^k,\ol{z}^k)\right)
	\\&
	-\alpha_k^{-2}\left(Q(\ol{x}^{k+1},y,z) - Q(x,\ol{y}^{k+1},\ol{z}^{k+1})\right),
\end{align*}
where \annotate.

Next, we sum these inequalities for $k = 0,\ldots,K-1$ and obtain the following:
\begin{align*}
	\mind{1em}
	\tau_x\alpha_{K}^{-2}\sqn{x^{K} - x}
	+\frac{1}{2\eta_y}\sqn{y^{K} - y}
	+\frac{1}{2\eta_z}\sqn{\hat{z}^{K} - z}
	+\frac{1}{2\eta_z}\sqn{\eta_z^{K} m^{K}}_{\mP}
	\\&\leq
	\tau_x\alpha_0^{-2}\sqn{x^0 - x}
	+\frac{1}{2\eta_y}\sqn{y^{0} - y}
	+\frac{1}{2\eta_z}\sqn{\hat{z}^{0} - z}
	+\frac{1}{2\eta_z}\sqn{\eta_z^0 m^0}_{\mP}
	+\frac{2n M^2 K}{\tau_xT}
	\\&
	-\frac{2\eta_y(K+2)^2}{9}\sqn{\tilde{x}^{K} - x^{K-1}}
	-\tfrac{1}{3}(K+2)\<x^{K-1} - \tilde{x}^K, y^{K} - y>
	\\&
	+\sum_{k=0}^{K-1}(\alpha_k^{-2}-\alpha_k^{-1})\left(Q(\ol{x}^{k},y,z) - Q(x,\ol{y}^k,\ol{z}^k)\right)
	\\&
	-\sum_{k=0}^{K-1}\alpha_k^{-2}\left(Q(\ol{x}^{k+1},y,z) - Q(x,\ol{y}^{k+1},\ol{z}^{k+1})\right)
	\\&\aleq{uses the Cauchy-Schwarz inequality and Young's inequality}
	\tau_x\alpha_0^{-2}\sqn{x^0 - x}
	+\frac{1}{2\eta_y}\sqn{y^{0} - y}
	+\frac{1}{2\eta_z}\sqn{\hat{z}^{0} - z}
	+\frac{1}{2\eta_z}\sqn{\eta_z^0 m^0}_{\mP}
	+\frac{2n M^2 K}{\tau_xT}
	\\&
	-\frac{2\eta_y(K+2)^2}{9}\sqn{\tilde{x}^{K} - x^{K-1}}
	+\frac{\eta_y (K+2)^2}{9}\sqn{x^{K-1} - \tilde{x}^K}
	+\frac{1}{4\eta_y}\sqn{y^K - y}
	\\&
	+\sum_{k=0}^{K-1}(\alpha_k^{-2}-\alpha_k^{-1})\left(Q(\ol{x}^{k},y,z) - Q(x,\ol{y}^k,\ol{z}^k)\right)
	\\&
	-\sum_{k=0}^{K-1}\alpha_k^{-2}\left(Q(\ol{x}^{k+1},y,z) - Q(x,\ol{y}^{k+1},\ol{z}^{k+1})\right)
	\\&\aeq{uses the fact that $\alpha_0 = 1$, which follows from \cref{eq:alpha_gamma_k}}
	\tau_x\sqn{x^0 - x}
	+\frac{1}{2\eta_y}\sqn{y^{0} - y}
	+\frac{1}{2\eta_z}\sqn{\hat{z}^{0} - z}
	+\frac{1}{2\eta_z}\sqn{\eta_z^0 m^0}_{\mP}
	+\frac{2n M^2 K}{\tau_xT}
	\\&
	-\frac{\eta_y(K+2)^2}{9}\sqn{\tilde{x}^{K} - x^{K-1}}
	+\frac{1}{4\eta_y}\sqn{y^K - y}
	-\alpha_{K-1}^{-2}\left(Q(\ol{x}^{K},y,z) - Q(x,\ol{y}^{K},\ol{z}^{K})\right)
	\\&
	+\sum_{k=1}^{K-1}(\alpha_k^{-2}-\alpha_k^{-1} - \alpha_{k-1}^{-2})\left(Q(\ol{x}^{k},y,z) - Q(x,\ol{y}^k,\ol{z}^k)\right)
	\\&\aeq{uses the definition of $\lambda_k$ in \cref{eq:lambda_k}}
	\tau_x\sqn{x^0 - x}
	+\frac{1}{2\eta_y}\sqn{y^{0} - y}
	+\frac{1}{2\eta_z}\sqn{\hat{z}^{0} - z}
	+\frac{1}{2\eta_z}\sqn{\eta_z^0 m^0}_{\mP}
	+\frac{2n M^2 K}{\tau_xT}
	\\&
	-\frac{\eta_y(K+2)^2}{9}\sqn{\tilde{x}^{K} - x^{K-1}}
	+\frac{1}{4\eta_y}\sqn{y^K - y}
	-\sum_{k=1}^{K}\lambda_k\left(Q(\ol{x}^{k},y,z) - Q(x,\ol{y}^k,\ol{z}^k)\right),
\end{align*}
where \annotate.
Further, we obtain the following:
\begin{align*}
	\mind{1em}
	\tau_x\alpha_{K}^{-2}\sqn{x^{K} - x}
	+\frac{1}{2\eta_y}\sqn{y^{K} - y}
	+\frac{1}{2\eta_z}\sqn{\hat{z}^{K} - z}
	+\frac{1}{2\eta_z}\sqn{\eta_z^{K} m^{K}}_{\mP}
	\\&\aleq{uses the convexity of $Q(x,y,z)$ in $x$ (follows from \Cref{ass:cvx}) and the concavity of $Q(x,y,z)$ in $(y,z)$, \cref{line:avg} of \Cref{alg}, and the fact that $\lambda_k\geq0$, which follows from \cref{eq:lambda_k,eq:alpha_gamma_k}}
	\tau_x\sqn{x^0 - x}
	+\frac{1}{2\eta_y}\sqn{y^{0} - y}
	+\frac{1}{2\eta_z}\sqn{\hat{z}^{0} - z}
	+\frac{1}{2\eta_z}\sqn{\eta_z^0 m^0}_{\mP}
	+\frac{2n M^2 K}{\tau_xT}
	\\&
	-\frac{\eta_y(K+2)^2}{9}\sqn{\tilde{x}^{K} - x^{K-1}}
	+\frac{1}{4\eta_y}\sqn{y^K - y}
	-\sum_{k=1}^{K}\lambda_k\left(Q(x_a^{K},y,z) - Q(x,y_a^K,z_a^K)\right)
	\\&\aeq{use the definition of $\lambda_k$ in \cref{eq:lambda_k} and the fact that $\alpha_0 = 1$, which follows from \cref{eq:lambda_k}}
	\tau_x\sqn{x^0 - x}
	+\frac{1}{2\eta_y}\sqn{y^{0} - y}
	+\frac{1}{2\eta_z}\sqn{\hat{z}^{0} - z}
	+\frac{1}{2\eta_z}\sqn{\eta_z^0 m^0}_{\mP}
	+\frac{2n M^2 K}{\tau_xT}
	\\&
	-\frac{\eta_y(K+2)^2}{9}\sqn{\tilde{x}^{K} - x^{K-1}}
	+\frac{1}{4\eta_y}\sqn{y^K - y}
	-\sum_{k=0}^{K-1}\alpha_k^{-1}\left(Q(x_a^{K},y,z) - Q(x,y_a^K,z_a^K)\right),
\end{align*}
where \annotate.
Next, we do rearranging and use \cref{eq:tau_x_eta_yz,eq:input}, which gives the following:
\begin{align*}
	\left(\sum_{k=0}^{K-1}\alpha_k^{-1}\right)\left(Q(x_a^{K},y,z) - Q(x,y_a^K,z_a^K)\right)
	\leq
	\frac{r}{3}\sqn{x}
	+\frac{6}{r}\sqn{y}
	+\frac{15\chi^2}{r}\sqn{z}
	+\frac{12n M^2 K}{r T}.
\end{align*}
Next, we divide both sides of the inequality by $\sum_{k=0}^{K-1}\alpha_k^{-1}$, which gives the following:
\begin{align*}
	Q(x_a^{K},y,z) - Q(x,y_a^K,z_a^K) \leq
	\left(\sum_{k=0}^{K-1}\alpha_k^{-1}\right)^{-1}
	\left(\frac{r}{3}\sqn{x}
	+\frac{6}{r}\sqn{y}
	+\frac{15\chi^2}{r}\sqn{z}
	+\frac{12n M^2 K}{r T}\right).
\end{align*}
Further, we can estimate $\sum_{k=0}^{K-1}\alpha_k^{-1}$ as follows:
\begin{align*}
	\sum_{k=0}^{K-1}\alpha_k^{-1}
	 & \aeq{uses \cref{eq:alpha_gamma_k}}
	\sum_{k=0}^{K-1} \frac{k+3}{3}
	=
	K + \frac{1}{3}\sum_{k=0}^{K-1} K
	=
	K + \frac{K(K-1)}{6}
	=
	\frac{K(K+5)}{6}
	\geq
	\frac{K^2}{6},
\end{align*}
where \annotate. Plugging this into the previous inequality gives
\begin{align*}
	Q(x_a^{K},y,z) - Q(x,y_a^K,z_a^K) & \leq
	\frac{6}{K^2}
	\left(\frac{r}{3}\sqn{x}
	+\frac{6}{r}\sqn{y}
	+\frac{15\chi^2}{r}\sqn{z}
	+\frac{12n M^2 K}{r T}\right)
	\\&=
	\frac{1}{K^2}
	\left(2r\sqn{x}
	+\frac{36}{r}\sqn{y}
	+\frac{90\chi^2}{r}\sqn{z}
	\right)
	+\frac{72n M^2}{r KT}
\end{align*}
which concludes the proof.\qed

\end{document}